\documentclass{article}

\usepackage{amssymb,amsmath,eucal}

\begin{document}

\def\R {{\mathbb R }}
\def\C {{\mathbb C }}
\def\K{{\mathbb K}}
\def\H{{\mathbb H}}
\def\Z{{\mathbb Z}}

\def\bm{{\mathbf m}}
\def\bl{{\mathbf l}}

\def\SS{{\cal S}_{p,q}}
\def\cS{{\cal S}}
\def\cB{{\cal B}}
\def\cD{{\cal D}}
\def\cR{{\cal R}}
\def\cH{{\cal H}}

\def\fD{{\frak D}}
\def\frk{\mathfrak k}

\def\wt{\widetilde}

\def\B{{\rm B}}

\def\OO{{\rm O}}
\def\SO{{\rm SO}}
\def\GL{{\rm GL}}
\def\SL{{\rm SL}}
\def\SU{{\rm SU}}
\def\Sp{{\rm Sp}}
\def\SOS{\SO^*}
\def\U{{\rm U}}

\def\Gr{{\rm Gr}}
\def\RP{{\R P}}
\def\CP{{\C P}}

\def\cS{{\cal S}}
\def\cE{{cal E}}

\def\fS{{\frak S}}

\def\b{C}
\def\T{T}

\def\ov{\overline}
\def\phi{\varphi}
\def\epsilon{\varepsilon}
\def\kappa{\varkappa}
\def\le{\leqslant}
\def\ge{\geqslant}

\def\konets{\hfill$\boxtimes$}

\newcommand{\Ker}{\mathop{\rm Ker}\nolimits}
  \renewcommand{\Im}{\mathop{\rm Im}\nolimits}
    \newcommand{\Dom}{\mathop{\rm Dom}\nolimits}
     \newcommand{\Indef}{\mathop{\rm Indef}\nolimits}
\newcommand{\graph}{\mathop{\rm graph}\nolimits}
\newcommand{\Mat}{\mathop{\rm Mat}\nolimits}

\def\const{{\rm const}}

\def\tto{\rightrightarrows}

\renewcommand{\Re}{\mathop{\rm Re}\nolimits}
\renewcommand{\Im}{\mathop{\rm Im}\nolimits}

\newcounter{sec}
 \renewcommand{\theequation}{\arabic{sec}.\arabic{equation}}

%\begin{figure}[t]

\def\kartinka{

\pagebreak

\unitlength=1mm

\begin{picture}(100,100)
%vertical lines
\put(10,50){\line(0,-1){43}}
\put(20,40){\line(0,-1){33}}
\put(30,30){\line(0,-1){23}}
\put(40,20){\line(0,-1){13}}
\put(50,10){\line(0,-1){3}}

%horizontal lines
\put(10,50){\line(-1,0){3}}
\put(20,40){\line(-1,0){13}}
\put(30,30){\line(-1,0){23}}
\put(40,20){\line(-1,0){33}}
\put(50,10){\line(-1,0){43}}

%circles
\multiput(10,80)(10,-10){8}{\circle*{1}}
\multiput(20,80)(10,-10){8}{\circle*{1}}
\multiput(30,80)(10,-10){8}{\circle*{1}}
\multiput(40,80)(10,-10){8}{\circle*{1}}

%axes
\put(60,60){\vector(0,1){30}}
\put(60,60){\vector(1,0){30}}
\put(90,65){$\sigma$}
\put(65,90){$\tau$}

 %squares
{\linethickness{1mm}
\multiput(10,50)(10,-10){5}%{\line(0,1){10}}
{\line(0,1){10}}
\multiput(20,50)(10,-10){5}%
{\line(0,1){10}}
\multiput(10,50)(10,-10) {5}%
{\line(1,0){10}}
\multiput(10,60)(10,-10) {5}%
{\line(1,0){10}}
\put(10,60){\line(0,1){10}}
\put(60,10){\line(1,0){10}}
\put(60,10){\line(0,-1){3}}
\put(70,10){\line(0,-1){3}}
\put(10,60){\line(-1,0){3}}
\put(10,70){\line(-1,0){3}}
}

%points
\multiput(11,51)(10,-10){5}
{%
\multiput(0,0)(1,0){9}%
{\multiput(0,0)(0,1){9}%
{\circle*{0.2}}%
}%
}

%principal series
\put(7,63){\line(1,-1){56}}

%labels
\put(70,60){\line(0,1){2}}
\put(70,63){1}
\put(60,70){\line(1,0){2}}
\put(63,70){1}

%center
\put(35,35){\circle*{2}}

%sdvig
\put(100,30){\vector(1,-1){10}}
\put(105,30){Shift}

%dots
\multiput(0,10)(0,10){8}\dots
\multiput(10,0)(10,0){11}\vdots

\end{picture}

%\end{figure}

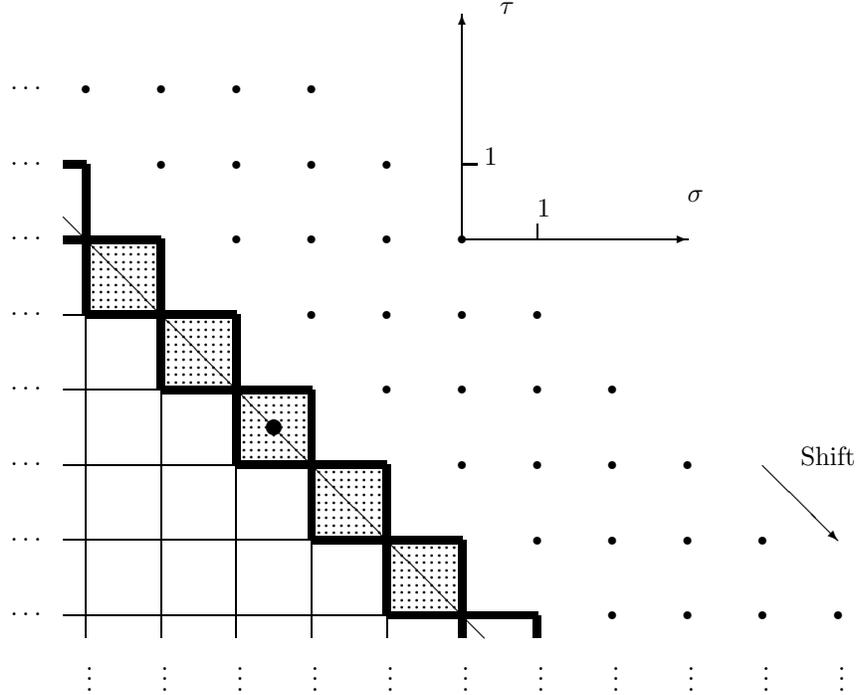
{\sc Figure.} {\sf Groups $\U(n)$. The case $n=5$.

 1. The dotted squares correspond to unitary representations
$\rho_{\sigma,\tau}$.

2. Vertical and horizontal rays in the south-west of Figure
correspond to nondegenerated
 highest weight and lowest weight representations.
Fat points  correspond to degenerated highest and lowest weight,
  representations, and also to unipotent representations. 
  The point
$(\sigma,\tau)=(0,0)$ corresponds
to the trivial one-dimensional representation.

3. In points of the thick segments,
we have some exotic unitary sub-factors.

4. The shift
$(\sigma,\tau)\mapsto(\sigma+1,\tau-1)$
 transform a projective representation
$\rho_{\sigma,\tau}$
to an equivalent representation.

5. The permutation of the axes
$(\tau,\sigma)\mapsto(\sigma,\tau)$
 gives a complex conjugate representation.

6. The symmetry with respect to the point $(-n/2,-n/2)$
(black circle) gives a dual representation
(for odd $n$ this point is a center of a dotted square;
for even $n$ this point is  a common vertex of two
dotted squares).

7. For $\sigma+\tau=n$ (the skew diagonal line)
our Hermitian
form is the standard $L^2$-product.
The represenations  $\rho_{\sigma,\tau}$ in this case are
in the untary principal series.

8. Linear (non-projective) reprsentations
of $\U(n,n)$ correspond to  the
family of parallel lines
$\sigma-\tau\in\Z$.}
\pagebreak
}
%\end{figure}

\begin{center}

{\Large \bf Notes on\\
 Sobolev spaces on compact classical groups
 \\ and Stein--Sahi representations}

%\vspace{22pt}

{\large\sc

Yurii A. Neretin}

\end{center}

{\small We discuss kernels
on complact classical groups
$G=\U(n)$, $\OO(2n)$, $\Sp(n)$
 defined by the formula
$K(z,u)=|\det(1-zu^*)|^s$. We obtain the explicit
Plancherel formula for these kernels and the
interval of positive-definiteness.
We also obtain explicit models for Sahi's
'unipotent' representations.}

\medskip

In \cite{Ner-determinant}, the author
proposed a class of natural kernels
on pseudo-Riemannian symmetric spaces
$G/H$
and conjectured that these kernels admit explicit
Plancherel formula.
The purpose of these notes
is to verify that at least in the case then
$G/H$ is a compact group 
$G\times G/G\sim\U(n)$, $\OO(2n)$, $\Sp(n)$
this Plancherel formula can be written.%
\footnote{Of course, the compact groups are the most
simple objects among symmetric spaces.
Nevertheless, our calculation (see also \cite{Ner-stein}) is
an important heurictic argument for 
a computability of the Plancherel measure in the general case.
 Indeed, Oshima \cite{Osh} recentely observed that
 difference equations for $c$-function
 are the same for different real forms of
 one complex symmetric space $G_\C/ H_\C$
 (and hence $c$-functions of different real forms differs
 by trigonometric factors;
 actual evaluation of this trigonometric factor is
 a highly nontrivial result of Oshima). It is natural
 to believe that this phenomenon survives in our situation.}

\smallskip

{\bf 0.1. A preliminary example.}
Recall the  construction   of unitary representations
of the complementary series of the group $\SL(2,\R)$
(Bargmann  \cite{Bar}, 1947, see also \cite{GN}).
Consider the unit circle $S^1$: $z=e^{i\phi}$
and the following Hermitian form in the space
of smooth functions  on $S^1$
$$
\langle f,g\rangle_s=
\int_0^{2\pi}\int_0^{2\pi}
\bigl|\sin (\phi/2)\bigr|^{1-s}
 f(\phi_1)\,\ov{g(\phi_2)}\,d\phi_1\,d\phi_2
.$$
This Hermitian form is invariant with respect
to the M\"obius transformations
 of the circle in the following sence.
Consider the group $\SL(2,\R)\simeq\SU(1,1)$
consisting of all complex $2\times 2$-matrices
$g=\begin{pmatrix}a&b\\ \ov b&\ov a\end{pmatrix}$,
satisfying $|a|^2-|b|^2=1$.
Then
the operators
$$
\rho_s\begin{pmatrix}a&b\\ \ov b&\ov a\end{pmatrix}
f(z)=f\bigl((\ov a+\ov b z)^{-1}(a+bz)\bigr) |\ov a+\ov b z|^{-1+s}
$$
%where
%$g=\begin{pmatrix}a&b\\ \ov b&\ov a\end{pmatrix}$,
%$|a|^2-|b|^2=1$; such matrices form
%the  group $\SL(2,\R)\simeq\SU(1,1)$,
%The operators $\rho_s$
preserve the form
$\langle f,g\rangle_s$, and moreover this property
uniquely determines our form.

For $-1<s<1$ our Hermitian form is positive definite.
This fact  follows from the
identity
$$
\langle\sum a_n e^{in\phi},b_n e^{in\phi}\rangle
=\const\sum \frac{\Gamma(n+(1-s)/2)}
                 {\Gamma(n+(1+s)/2)}
          a_n\ov b_n
.$$
We observe, that for $-1<s<1$,
all the gamma-coefficients are positive,
and hence $\langle f,f\rangle_s$ is positive
for all $f$.

It is natural to consider  the completion $H_s$
 of the space $C^\infty(S^1)$
with respect to this inner product.   Obviously,
\begin{multline*}
f(\phi)=\sum a_n e^{in\phi}\in H_s
\quad\Longleftrightarrow\quad
\sum
\frac{\Gamma(n+(1-s)/2)}
     {\Gamma(n+(1+s)/2)}   |a_n|^2<\infty
\quad\Longleftrightarrow
\\ \Longleftrightarrow   \quad
 \sum \bigl(1+|n|\bigr)^{-s}|a_n|^2<\infty
.\end{multline*}
The last condition shows that $H_s$ is a Sobolev space.

This construction is the simplest example
of a representation that is unitary by
some nontrivial reason.
In fact, Hilbert spaces with various inner product
having the form
$$
\langle f,g\rangle
  =\iint_{X\times X} K(x,y)  f(x)\,\ov g(y)\,dx\,dy
,$$
where $K(x,y)$ is a distribution,
are usual in the representation theory.

For instance (see Vilenkin, \cite{Vil}, X.2),
conformally invariant inner products in spaces
of functions on a spere $|x_1|^2+\dots+|x_n|^2=1$
are given by the formula
\begin{equation}
\langle f,g\rangle=
\iint_{S^{n-1}\times S^{n-1}}
\|x-y\|^{-\lambda} f(x)\ov {g(y)}\,dx\,dy
\end{equation}
and this inner product is positive definite
iff $0<\lambda<n-1$.

The  Hilbert spaces defined by the inner products
(0.1)
 are Sobolev spaces.
In the representation theory also there are few cases
(related to the groups $\U(1,n)$, $\Sp(1,n)$)
when some anisotropic Sobolev spaces arise
in a natural way,
but usually the situation is more complicated
(some discussion of functional-theoretical
 phenomena is contained in \cite{NO}).

\smallskip

{\bf 0.2. Stein kernels.}
Consider the space $\Mat_n$ consisting
of $n\times n$ complex matrices. Consider
the following Hermitian form on $C^\infty(\Mat_n)$
\begin{equation}
\langle f,g\rangle=
\iint_{\Mat_n\times\Mat_n}
|\det(x-y)|^{-2n+\sigma} f(x)\ov{g(y)}  \,dx\,dy
\end{equation}
(for $\sigma>2n$ the integral is well defined,
further we can consider
the meromorphic continuation in  $\sigma$).
The Hermitian form (0.2)
 is invariant with respect to the transformations
$$
\rho_\sigma\begin{pmatrix}a&b\\c&d\end{pmatrix}
f(x)=
f((a+xc)^{-1}(b+xd))|\det(a+xc)|^{-2n-\sigma}
,$$
where
$\begin{pmatrix}a&b\\c&d\end{pmatrix}\in\GL(2n,\C)$
 is an arbitrary
invertible $(n+n)\times(n+n)$ matrix.

\smallskip

{\sc Stein's theorem.} {\it The inner product (0.2)
is positive definite iff $-1<\sigma<1$.}

\smallskip

Vogan \cite{Vog}, 1986, extended Stein's construction
to the groups $\GL(n)$ over real numbers $\R$
and quaternions $\H$, see also \cite{Sah1}.
Sahi \cite{Sah2},  \cite{Sah3}
extended the construction to other series
of classical groups, precisely to
the groups
$\OO(2n,2n)$, $\U(n,n)$, $\Sp(n,n)$,
$\Sp(2n,\R)$, $\SOS(4n)$, $\Sp(4n,\C)$, and
$\OO(2n,\C)$, 

\smallskip

The inner products (0.1) and (0.2) seem similar, but
the kernel $\|x-y\|^{-\lambda}$ has a singularity
on the diagonal $x=y$; on the contrary
$|\det(x-y)|^{-2n-\sigma}$ has a singularity
on  a  complicated surface in $\Mat_n\times \Mat_n$
containing the diagonal.
In particular, the Hilbert space  defined by
the Stein inner product
(0.2) is  not a Sobolev space
in the standard sence.

Similarily, various natural integral operators
that appear in the representation theory usually
are not pseudodifferential operators
 in the standard sence.

\smallskip

{\bf 0.3. Sobolev kernels on the orthogonal
groups $\OO(2n)$.}
Consider the orthogonal group
$\OO(2n)$, i.e., the group of real $2n\times 2n$-matrices
$h$ satisfying $h^t h=1$.
We consider the Hermitian form on $C^\infty(\OO(2n))$
given by
\begin{equation}
\langle F_1, F_2\rangle_\lambda=
\iint_{\OO(2n)\times\OO(2n)}
|\det(u-v)|^{\lambda} F_1(u)\ov{F_2(v)}\,d\mu(u)\,d\mu(v)
,\end{equation}
where $\mu$ is the Haar measure on $\OO(2n)$.
This form is well defined for $\lambda>0$, for $\lambda<0$
we consider its meromorphic continuation.

\smallskip

In this paper, we obtain the explicit expression
of this form in characters (Theorem 3.2)
 and also show that for
 $-n<\lambda<-n+1$ our form is positive definite

\smallskip

We consider the group $\OO(2n,2n)$ consisting
of $(2n+2n)\times (2n+2n)$-matrices
$\begin{pmatrix}a&b\\c&d\end{pmatrix}$
satisfying the condition
$$
\begin{pmatrix}a&b\\c&d\end{pmatrix}^t
\begin{pmatrix}1&0\\0&-1\end{pmatrix}
\begin{pmatrix}a&b\\c&d\end{pmatrix}
=
\begin{pmatrix}1&0\\0&-1\end{pmatrix}
,$$
where the symbol ${}^t$ denotes the transposition of a matrix.

For an orthogonal matrix $h\in\OO(2n)$ and
$\begin{pmatrix}a&b\\c&d\end{pmatrix}
\in\OO(2n,2n)$, we have
$$(a+hc)^{-1}(b+hd)\in\OO(2n)
.$$
 It is easy to verify that
the Hermitian form (0.3) is invariant with respect
to the operators $C^\infty(\OO(2n))\to C^\infty(\OO(2n))$
given by
$$\rho_\lambda
\begin{pmatrix}a&b\\c&d\end{pmatrix}
F(h)=F\bigl((a+hc)^{-1}(b+hd)\bigr)
\det(a+hc)^{-2n+1-\lambda}
.$$

Thus our Theorem 3.2 gives another proof of Sahi's
theorem about existence of Stein-type series for
the groups $\OO(2n,2n)$.
Also, we obtain new models of the "unipotent representations"
of the groups $\OO(2n,2n)$ 
(see Sahi, \cite{Sah4}, Dvorsky and Sahi,
\cite{DS1}, \cite{DS2}).

\smallskip

{\bf 0.4. Structure of the paper.}
We consider the groups $\U(n)$, $\OO(2n)$, $\Sp(n)$
in Sections 2, 3, 4 respectively.
We also obtain models of Stein--Sahi representations
for the groups $\U(n,n)$, $\OO(2n,2n)$, $\Sp(n,n)$
In particular,
we obtain an independent proof
 of the corresponding Sahi's results.

 In all the cases,
we reduce the problem to evaluation of some determinants,
necessary determinant calculations are collected
in  preliminary Section 1 (our calculations are more
uniform than it seems at the first glance, see \cite{Ner-stein}.

In Section 5, we discuss models of 'unipotent representations'
from \cite{Sah4}, \cite{DS1}, \cite{DS2}.

\smallskip

{\bf 0.5. Notation.}
We denote the Haar measure on $\U(n)$, $\SO(n)$, $\Sp(n)$
 by $\mu$. We assume that the measure of the whole group
is 1.

The Pochhammer symbol  is given by
\begin{equation}
(a)_n=a(a+1)\dots(a+n-1)=\frac{\Gamma(a+n)}{\Gamma(a)}
.\end{equation}
The second expression also has sence for negative $n$.

\smallskip

{\bf Acknowledgements.}
I am grateful to A.Dvorsky for his explanations on 
'unipotent representations'.
 This work was done during my visit
to the Schr\"odinger Institute, Vienna, winter 2001--2002.
I thank the administrators of the Institute for
their hospitality.

\vspace{22pt}

{\large\bf 1. Some determinant identities}

\medskip

\stepcounter{sec}
\setcounter{equation}{0}

\nopagebreak

\def\palka{\vphantom{1^G}}

Let $A=\{a_{kl}\}$ be a square $n\times n$ matrix,
$k,l=1,2,\dots,n$,
We denote its determinant by
$\det_{k,l} a_{kl}$.

\smallskip

{\bf 1.1. A determinant of Cauchy type.}
 Recall that the {\it  Cauchy determinant}
(see, for instance, \cite{Kra})
is given by
\begin{equation}
\det\limits_{kl} \frac 1{x_k+y_l}=
\frac{
\prod_{1\le k <l\le n}(x_k-x_l)
\cdot
\prod_{1\le k <l\le n}(y_k-y_l)}
{\prod\limits_{1\le k, l\le n} (x_k+y_l)}
.\end{equation}

The following variant of the Cauchy determinant,
 is also well known.

\smallskip

{\sc Lemma 1.1.}
\begin{multline}
\det\begin{pmatrix}
1&1&1& \dots & 1 \\
\frac 1{x_1+b_1} & \frac 1{x_2+b_1} & \frac 1{x_3+b_1}&
    \dots & \frac {1\vphantom{1^G}}
{x_n+b_1}\\
\frac {1\vphantom{1^G}}{x_1+b_2} & \frac 1{x_2+b_2} & \frac 1{x_3+b_2}&
    \dots & \frac 1{x_n+b_2}\\
\vdots&\vdots&\vdots&\ddots &\vdots\\
\frac {1\vphantom{1^G}}{x_1+b_{n-1}} & \frac 1{x_2+b_{n-1}} & \frac 1{x_3+b_{n-1}}&
    \dots & \frac 1{x_n+b_{n-1}}\\          \end{pmatrix}
=\\=
\frac{
\prod_{1\le k <l\le n} (x_k-x_l)
 \prod_{1\le \alpha<\beta\le n-1} (b_\alpha-b_\beta)}
{\prod\limits_{\begin{smallmatrix} 1\le k \le n\\
                   1\le \alpha \le n-1 \end{smallmatrix}}
(x_k+b_\alpha)}
.\end{multline}

{\sc Proof.}  Let $\Delta$ be the Cauchy determinant (1.1).
Then
$$
y_1\Delta=\begin{pmatrix}
\frac {y_1}{x_1+y_1}&  \frac {y_1}{x_2+y_1}&\dots & \frac {y_1}{x_1+y_1}\\
\frac 1{x_1+y_2}&  \frac 1{x_2+y_2}&\dots & \frac 1{x_n+y_2}\\
\vdots&\vdots&\ddots &\vdots\\
\frac 1{x_1+y_n}&  \frac 1{x_2+y_n}&\dots & \frac 1{x_n+y_n}
\end{pmatrix}
.$$
We consider  $\lim\limits_{y_1\to \infty} y_1 \Delta$
and substitute $y_{\alpha+1}=b_\alpha$.

\smallskip

{\bf 1.2. One standard determinant.}
The following determinant is equivalent to
Lemma 3 from Krattenthaller, \cite{Kra}.

{\sc Lemma 1.2.}

\begin{multline}
\!\!\!\!\!\!\!\!\!
\det\begin{pmatrix}
1&1&1& \dots & 1 \\
\frac {x_1+b_1\palka} {x_1+a_1} & \frac {x_2+a_1}{x_2+b_1} &
     \frac {x_3+a_1}{x_3+b_1}&
    \dots & \frac {x_n+a_1}{x_n+b_1}\\
\frac {(x_1+a_1)\palka(x_1+a_2)}{(x_1+b_1)(x_1+b_2)}
  & \frac {(x_2+a_1)(x_2+a_2)}{(x_2+b_1)(x_2+b_2)}
  & \frac {(x_3+a_1)(x_3+a_2)}{(x_3+b_1)(x_3+b_2)}&
    \dots & \frac {(x_n+a_1)(x_n+a_2)}{(x_n+b_1)(x_n+b_2)}\\
\vdots & \vdots & \vdots &\ddots &\vdots\\
\frac {\prod\limits_{1\le m\le n-1} (x_1+a_m)}
{\prod\limits_{1\le m\le n-1} (x_1+b_m)}&
\frac {\prod\limits_{1\le m\le n-1} (x_2+a_m)}
     {\prod\limits_{1\le m\le n-1} (x_2+b_m)}&
\frac {\prod\limits_{1\le m\le n-1} (x_3+a_m)}
   {\prod\limits_{1\le m\le n-1} (x_3+b_m)} &
\dots&
\frac {\prod\limits_{m:\,1\le m\le n-1} (x_n+a_m)}
 {\prod\limits_{m:\,1\le m\le n-1} (x_n+b_m)}
\end{pmatrix}
=\\=
\frac{\prod\limits_{1\le k<l\le n} (x_k-x_l)
  \prod\limits_{1\le \alpha\le\beta\le n-1} (a_\alpha-b_\beta)}
{\prod\limits_{1\le k\le n, 1\le\beta\le n-1} (x_k+b_\beta)}
.\end{multline}

{\sc Proof.}
Decomposing a matrix element into the sum
of partial fractions, we obtain
\begin{multline*}
\frac{(x_k+a_1)\dots (x_k+a_\alpha)}
     {(x_k+b_1)\dots (x_k+b_\alpha)}=
1+\sum\limits_{1\le \beta \le \alpha}
\frac{(a_1-b_\beta)(a_2-b_\beta)\dots(a_\alpha-b_\beta)}
    {(b_1-b_\beta)\dots (b_{\beta-1}-b_\beta)
            (b_{\beta-1}-b_\beta)}
\cdot \frac 1{x_k+b_\beta}
\end{multline*}

%=\\=
%\frac{(a_1-b_\alpha)(a_2-b_\alpha)\dots(a_\alpha-b_\alpha)}
%  {(b_1-b_\alpha)(b_2-b_\alpha)\dots (b_{\alpha-1}-b_\alpha)}
%\cdot \frac 1{x_k+b_\alpha}
%+ \dots
%\end{multline*}

We observe that the  $(\alpha+1)$-th row
is a linear combination of the rows
$$
\begin{matrix}
\Bigl(&1&1&\dots&1&\Bigr), \\
\Bigl(&\frac 1{x_1+b_1}& \frac 1{x_2+b_1} &\dots&\frac 1{x_n+b_1} &\Bigr),\\
\dots&\dots&\dots&\dots&\dots&\dots\\
\Bigl(&\frac 1{x_1+b_\alpha}&
 \frac 1{x_2+b_\alpha} &\dots&\frac 1{x_n+b_\alpha} &\Bigr).\\
\end{matrix}
.$$
Thus our determinant is
$$%\begin{multline*}
\prod_{\alpha=1}^{l-1}
\frac{\prod_{j=1}^{\alpha}(a_j-b_\alpha)}
{\prod_{j=1}^{\alpha-1}(b_j-b_\alpha)}
\cdot
\det
\begin{pmatrix}
1&1&\dots&1\\
\frac 1{x_1+b_1}& \frac 1{x_2+b_1} &\dots&\frac 1{x_n+b_1} \\
\vdots&\vdots&\ddots&\vdots\\
\frac 1{x_1+b_\alpha}&
 \frac 1{x_2+b_\alpha} &\dots&\frac 1{x_n+b_\alpha}\\
\end{pmatrix}
.$$%\end{multline*}
and we reduce the evaluation
of  our determinant to Lemma 1.1.

 \smallskip

{\bf 1.3. Two determinants.}

\smallskip

{\sc Lemma 1.3.} {\it Concider the
$n\times n$ matrix $Q$ having elements
$$
q_{lk}=\frac{(x_k+a_1)(x_k+a_2)\dots (x_k+a_{l})}
            {(x_k+b_1)(x_k+b_2)\dots (x_k+b_{l})}
                           -
      \frac{(x_k-a_1)(x_k-a_2)\dots (x_k-a_{l})}
     {(x_k-b_1)(x_k-b_2)\dots (x_k-b_{l})}
.$$
Then              }
$$
\det Q=\frac
  {2^n \prod\limits_{1\le k\le n} x_k \cdot
      \prod\limits_{1\le k<l \le n} (x_k^2-x_l^2)
  \cdot \prod\limits_{1\le\alpha <\beta\le n} (b_\alpha+b_\beta)
   \cdot \prod\limits_{1\le \alpha \le \beta\le n} (a_\alpha-b_\beta)}
  {\prod\limits_{1\le k \le n, 1\le \alpha \le n-1} (x_k^2-b_\alpha^2)}
.$$

{\sc Proof.} We expand a matrix element in a sum of
partial fractions.
$$q_{kl}=
\sum_{\alpha=1}^l
\frac{ \prod_{j=1}^l (a_j-b_\alpha)}
{\prod_{j=1}^{\alpha-1}(b_j-b_\alpha)}
\Bigl\{
\frac {1}{x_k+b_\alpha}+\frac {1}{x_k-b_\alpha}
\Bigr\}
.$$
We write
$$
\frac {1}{x_k+b_\alpha}+\frac {1}{x_k-b_\alpha}
         =\frac{2x_k}{x_k^2-b_\alpha^2}
$$
and observe that an $l$-th row is a linear combination
of vector-rows
$$
\begin{matrix}
\Bigl(&\frac{x_1}{x_1^2-b_1^2}&\dots &\frac{x_k}{x_n^2-b_1^2}& \Bigr),\\
\dots&\dots&\dots&\dots&\dots \\
\Bigl(&\frac{x_1}{x_1^2-b_l^2}&\dots &\frac{x_k}{x_n^2-b_l^2}& \Bigr).\\
\end{matrix}
$$
Thus our determinant is
$$
\prod_{l=1}^{n}
\frac{ \prod_{j=1}^l (a_j-b_\alpha)}
{
\prod_{\alpha=1}^{\alpha-1}(b_j-b_\alpha)
}
\cdot
\prod_{j=1}^n (2x_j)
\cdot
\det_{1\le j,l\le n} \frac 1{x_j^2-b_l^2}
.$$
and the last factor is reduced to the classical Cauchy determinant (1.1).

\smallskip

{\sc Lemma 1.4.}
\begin{multline}
\det\limits_{k,l}
\Biggl\{\frac{(x_k+a_1)(x_k+a_2)\dots (x_k+a_{l})}
            {(x_k+b_1)(x_k+b_2)\dots (x_k+b_{l})}
                           +
      \frac{(x_k-a_1)(x_k-a_2)\dots (x_k-a_{l})}
     {(x_k-b_1)(x_k-b_2)\dots (x_k-b_{l})}
\Biggr\}
=\\=
 \frac { 2
      \prod\limits_{1\le k<l \le n} (x_k^2-x_l^2)
  \cdot \prod\limits_{1\le\alpha \le\beta\le n-1} (b_\alpha+b_\beta)
   \cdot \prod\limits_{1\le \alpha \le \beta\le n-1} (a_\alpha-b_\beta)}
  {\prod\limits_{1\le k \le n, 1\le \alpha \le n-1} (x_k^2-b_\alpha^2)}
.\end{multline}

{\sc Remark.} In particular the first row of our matrix is
$\begin{pmatrix} 2&\dots&2\end{pmatrix}$.

\smallskip

{\sc Proof.}
Decomposing a matrix element  into a sum of prime fractions,
we obtain
$$q_{kl}=
2+\sum\limits_{\alpha=1}^l
\frac{\prod_{j=1}^{\alpha} (a_j-b_\alpha)}
    { \prod_{j=1}^{\alpha-1} (b_j-b_\alpha)}
     \Bigl(
     \frac 1{x_k+b_j} - \frac 1{x_k-b_j}\Bigr)
$$
Since,
$$  \frac 1{x_k+b_l} - \frac 1{x_k-b_l}   =
\frac {-2b_k}{x_k^2-b_l^2},
$$
the $l$-th row  is a linear combination of the vectors-rows
$$
\begin{matrix}
\Bigl(& 2&\dots&2&\Bigr)\\
\Bigl(&\frac{b_1}{x_1^2-b_1^2}&\dots &\frac{b_1}{x_n^2-b_1^2}& \Bigr)\\
\dots&\dots&\dots&\dots&\dots \\
\Bigl(&\frac{b_{l-1}}{x_1^2-b_{n-1}^2}&\dots
        &\frac{b_n}{x_n^2-b_{n-1}^2}& \Bigr)\\
\end{matrix}
$$
and  we obtain
\begin{multline*}
2\prod\limits_{1\le j \le {n-1}} 2b_j
\prod\limits_{1\le l\le n-1}
\frac{(a_1-b_l)(a_2-b_l)\dots (a_l-b_l)}
   {(b_1-b_l)(b_2-b_l) \dots (b_{l-1}-b_l)}
\times\\ \times
\det\limits_{k,l}
\begin{pmatrix}
1&1&\dots&1\\
\frac 1{x_1^2-b_1^2}& \frac 1{x_2^2-b_1^2}&\dots&\frac 1{x_n^2-b_1^2}\\
\vdots&\vdots&\ddots&\vdots\\
\frac 1{x_1^2-b_{n-1}^2}&
      \frac 1{x_2^2-b_{n-1}^2}&\dots&\frac 1{x_n^2-b_{n-1}^2}\\
\end{pmatrix}
.\end{multline*}
It remains to apply Lemma 1.1.

%\Bigl\{\frac 1{x_k^2-b_l^2}\Bigr\}

%Proof of the statement a) is the same,
%we only use Lemma 1.1 unstead of Cauchy determinant.
%\konets

\bigskip

{\bf\large 2. Sobolev kernel on the unitary group.}

\medskip

\stepcounter{sec}
\setcounter{equation}{0}

{\bf 2.1. Definition of kernel.}
Let $z$ be an $n\times n$ matrix
with norm $< 1$.
For $\sigma\in\C$,
we define the function $\det(1-z)^\sigma$ by
$$
\det(1-z)^\sigma:=\det\Bigl[ 1-\sigma z+\frac{\sigma(\sigma-1)}{2!} g^2-
  \frac{\sigma(\sigma-1)(\sigma-2)}{3!} g^3+\dots \Bigr]
.$$
We also define this function for $z$ satisfying $\|z\|=1$, $\det(1-z)\ne 0$
being
$$
\det(1-z)^\sigma:=\lim_{u\to z, \,|u\|<1}\det(1-u)^\sigma
.$$
The expression $\det(1-z)^\sigma$ is continuous
on the domain $\|z\|\le 1$ except the surface
$\det(1-z)=0$.

We denote by $\det(1-z)^{\{\sigma|\tau\}}$ the function
$$
\det(1-z)^{\{\sigma|\tau\}}:=
  \det(1-z)^\sigma \det(1-\ov z)^\tau
.$$
We define the function
$\ell_{\sigma,\tau}(g)$
on the unitary group $\U(n)$ by
$$
\ell_{\sigma,\tau}(g): =
2^{-(\sigma+\tau)n}
\det(1-z)^{\{\sigma|\tau\}}
.$$
Obviously,
\begin{equation}
\ell_{\sigma,\tau} (h^{-1}gh)=\ell_{\sigma,\tau}  (g)
\qquad\text{for $g$, $h\in\U(n)$}
.\end{equation}

{\sc Lemma 2.1.} {\it Let $e^{i\psi_1}$, \dots, $e^{i\psi_n}$
be the eigenvalues of $g\in\U(n)$; we assume $0\le \psi_k<2\pi$.
Then        }
$$
\ell(g)=
\exp\bigl\{(\sigma-\tau)\sum_k(\psi_k-\pi)/2\bigr\}
\prod_{k=1}^n \sin^{\sigma+\tau} \frac {\psi_k}2
.$$

{\sc Proof.} It is sufficient to verify this statement
for diagonal matrices, or equivalenly we can check
the identity
$$(1-e^{i\psi})^{\{\sigma|\tau\}}=
\exp\bigl\{(\sigma-\tau)(\psi-\pi)/2\bigr\}
\sin^{\sigma+\tau} \frac {\psi}2
.$$
We have
$$
\frac12(1-e^{i\psi})=
\exp\{i(\psi-\pi)/2\}\sin\frac \psi 2
.$$
Further, in the equality
$$
2^{-\sigma}(1-e^{i\psi})^\sigma=
\exp\{i\sigma(\psi-\pi)/2\}\sin^\sigma\frac \psi 2
,$$
the both sides are real-analytic on $(0,2\pi)$
and the substitution $\psi=\pi$ gives 1 in the both sides.
\hfill $\square$

\smallskip

We also define  the kernel
$L_{\sigma,\tau}(g,h)$ on $\U(n)$ by
\begin{equation}
L_{\sigma,\tau}(g,h)=\ell_{\sigma,\tau}(gh^{-1})
.\end{equation}
Obviously, this kernel is invariant with respect to
left and right shifts on $\U(n)$, i.e.,
$$L_{\sigma,\tau} (r_1g r_2, r_1 h r_2)=
  L_{\sigma,\tau} (g , h )
\qquad\text{for $g$, $h$, $r_1$, $r_2\in\U(n)$}
.$$

{\bf 2.2. Characters, } see Weyl book \cite{Wey}.
 The set of finite dimensional representations
of $\U(n)$ is parametrized by collections of integers
(signatures)
$$
\bm: \quad m_1> m_2 >\dots > m_n
.$$
The character $\chi_\bm$ of representation $\pi_\bm$
(a Schur function)
corresponding to a signature $\bm$
is given by
\begin{equation}
\chi_\bm(g)=
   \frac    {\det_{k,j=1,2,\dots,n} \bigl\{e^{i m_j \psi_k}\bigr\} }
            {\det_{k,j=1,2,\dots,n} \bigl\{e^{i(j-1) \psi_k}\bigr\}    }
,\end{equation}
where $e^{i\psi_k}$ is the eigenvalues of $g$.
  Recall that the denominator admits decomposition
\begin{equation}
\prod_{l<k} (e^{i\psi_l}-e^{i\psi_k} ).
\end{equation}

The dimension  of $\pi_m$ is
\begin{equation}
\dim \pi_\bm=\chi_\bm(1)=
\frac{\prod_{0\le\alpha<\beta\le n} (m_\alpha-m_\beta)}
{\prod_{j=1}^n j!}
.\end{equation}

A function $F(g)$ on $\U(n)$ is {\it central}
if it satisfies the identity
$F(h^{-1} g h)=F(g)$.

Consider the Haar measure $\mu$ on $\U(n)$ normalized by the condition:
the measure of the whole group is 1.
For a central function on $\U(n)$,
 the following {\it Weyl integration formula}
holds
\begin{multline}
\int\limits_{\U(n)} F(g) \,d\mu(g)
= \\ =
\frac 1{(2\pi)^n n!}
 \!\!\!\!\!\!\!\!
\int\limits_{0<\psi_1<2\pi}\dots\int\limits_{0<\psi_n<2\pi}
\!\!\!\!\!\!\!\!
F\bigl(
\mathrm{diag}(e^{i\psi_1},\dots,e^{i\psi_n})
\bigr)
\Bigl| \prod_{m<k} (e^{i\psi_m}-e^{i\psi_k} ) \Bigr|^2
\,\prod_{k=1}^n d\phi_k
,\end{multline}
where $\mathrm{diag}(\cdot)$
is a diagonal matrix with given entries.

A  central function $F\in L^2(\U(n))$
 admits an expansion in characters,
$$
F(g)=\sum_\bm c_\bm\chi_m
.$$
where the summation is given over all the signatures $\bm$
and  the coefficients
 $c_\bm$ are the $L^2$-inner products
$$
c_\bm=\int_{\U(n)} F(g) \overline{\chi_\bm(g)}\,d\mu(g)
.$$
Applying formula (2.6), explicit expression (2.3) for characters,
and formula (2.4) for the denominator, we obtain
\begin{multline}
c_{\bm}=
\frac 1{(2\pi)^n n!}
\int\limits_{0<\psi_1<2\pi}\dots\int\limits_{0<\psi_n<2\pi}
F\Bigl(\mathrm{diag}\bigl\{e^{i\psi_1},\dots,e^{i\psi_n}\bigr\}\Bigr)
\times\\ \times
 {\det_{k,j=1,2,\dots,n} \bigl\{e^{i(j-1) \psi_k}\bigr\}    }
      \det_{k,j=1,2,\dots,n} \bigl\{e^{-i m_j \psi_k}\bigr\}
\,\prod_{k=1}^n d\phi_k
.\end{multline}

For calculation of such expressions, we will use the following
evident Lemma  (see, for instance, \cite{Ner-determinant},)

\smallskip

{\sc Lemma 2.2.} {\it Let $X$ be a set,}
\begin{multline*}
\int_{X^n} \prod_{k=1}^n f(x_k) \,
\det\limits_{k,l=1,\dots n} \{u_l(x_k)\}
\det\limits_{k,l=1,\dots n} \{v_l(x_k)\}
\prod_{j=1}^n dx_j
=\\
= n!
\det\limits_{l,m=1,\dots, n}\Bigl\{
\int_X f(x)u_l(x)v_m(x)\,dx
\Bigr\}
\end{multline*}

{\bf 2.3. Lobachevsky beta-integrals.}
We will use two following integrals,
see \cite{GR}, 3.631,1,  3.631,8,
\begin{align}
\int_0^\pi \sin^{\mu-1} (\phi) \, e^{ibx}\,dx=
\frac{2^{1-\mu}\pi \Gamma(\mu) e^{ib\pi/2}}
    {\Gamma\bigl((\mu+b+1)/2\bigr)
     \Gamma\bigl((\mu-b+1)/2\bigr) }
\\
\int_0^\pi \sin^{\mu-1} (\phi) \, \cos(bx)\,dx=
\frac{2^{1-\mu}\pi \Gamma(\mu) \cos(b\pi/2)}
    {\Gamma\bigl((\mu+b+1)/2\bigr)
     \Gamma\bigl((\mu-b+1)/2\bigr) }
%\\
%\int_0^\pi \sin^{\mu-1} \phi\, \sin(bx)\,dx=
%\frac{2^{1-\mu}\pi \Gamma(\mu) \sin(b\pi/2)}
%    {\Gamma\bigl((\mu+b+1)/2\bigr)
%     \Gamma\bigl((\mu-b+1)/2\bigr) }
\end{align}

In some sence, our integral evaluations below,
2.4, 3.3, 4.3, are multivariate analogs
of these  integrals.

{\bf 2.4. Expansion of the function $\ell_{\sigma,\tau}$.}

\smallskip

{\sc Theorem 2.3.}  {\it Let $\Re (\sigma+\tau)<1$.
Then}
\begin{align}
\ell_{\sigma,\tau}(g)&= \nonumber
\\
&=
\frac{(-1)^{n(n-1)/2}\sin^n(\pi\sigma)
2^{-(\sigma+\tau)n}}{\pi^n}
\prod_{j=1}^n\Gamma(\sigma+\tau+j)
\times  \notag \\
&\qquad\times
\sum\limits_\bm\Biggl\{
\prod\limits_{1\le\alpha<\beta\le n} (m_\alpha-m_\beta)
\prod\limits_{j=1}^n\frac{\Gamma(-\sigma+m_j-n+1)}
    {\Gamma(\tau+m_j+1)}\chi_\bm(g)
\Biggr\}
=\\
&=
(-1)^{n(n-1)/2}
2^{-(\sigma+\tau)n}
\prod_{j=1}^n\Gamma(\sigma+\tau+j)
\times\notag
\\ &\quad\times
\sum\limits_\bm\Biggl\{
\frac{(-1)^{\sum m_j}
\prod\limits_{1\le\alpha<\beta\le n} (m_\alpha-m_\beta)}
{\prod\limits_{j=1}^n\Gamma(\sigma-m_j+n)
    \Gamma(\tau+m_j+1)}\chi_\bm(g) \Biggr\}
.\end{align}

{\sc Proof.}
We must evaluate the
inner product
$$
\int_{\U(n)} \ell_{\sigma,\tau}(g)\,\ov{\chi_\bm(g)}\,d\mu(g)
$$
Applying (2.7), we obtain
\begin{multline*}
\frac1{(2\pi)^n\,n!}
\int\limits_{0<\psi_k<2\pi}
\prod\limits_{j=1}^n
\Bigl[\sin^{\sigma+\tau}
  \bigl( \psi_j/ 2\bigr)\cdot
 \exp\bigl\{i(\sigma-\tau)(\psi_j-\pi)/2\bigr\}
\Bigr]
\times\\
\times \det\limits_{1\le k,l\le n}
    \{e^{-im_k\psi_l}\}
\cdot\det\limits_{1\le k,l\le n}
    \{e^{i(k-1)\phi_l}\}
\prod\limits_{l=1}^n d\phi_l
.\end{multline*}
By Lemma 2.2, we reduce this integral to
$$
\frac1{(2\pi)^n}
\det\limits_{1\le k,j\le n} I(k,j)
,$$
where
$$
I(k,j)=e^{-i(\sigma-\tau)\pi/2}
\int_0^{2\pi}
\sin^{\sigma+\tau}( \psi/ 2)
\cdot
 \exp\bigl\{i(\,(\sigma+\tau)/2+k-1-m_j)\bigr\}
\,d\phi
.$$
We apply (2.8)
and obtain
$$
I(k,j)=\frac{2^{1-\sigma-\tau}\pi\Gamma(\sigma+\tau+1)\, (-1)^{k-1-m_j} }
            {\Gamma(\sigma+k-m_j)\Gamma(\tau-k+m_j+2)}
$$
Applying standard formulae for $\Gamma$-function,
we obtain
\begin{multline*}
I(k,j)=2^{1-\sigma-\tau}\Gamma(\sigma+\tau+1)\, \sin(-\sigma\pi)
\cdot
           \frac{ \Gamma(-\sigma+m_j-k+1)}
            {\Gamma(\tau+m_j-k+2)}
                 =\\=
2^{1-\sigma-\tau}\Gamma(\sigma+\tau+1)\, \sin(-\sigma\pi)
\cdot
           \frac{ \Gamma(-\sigma+m_j-n+1)}
                  {\Gamma(\tau+m_j-n+2)}
\, \cdot
       \,  \boxed {  \frac{ (-\sigma+m_j-n+1)_{n-k}}
                  { (\tau+m_j-n+2)_{n-k}}
               }
\end{multline*}
The factors outside the box do not depend on
on $k$. Thus, we must evaluate the determinant
$$
\det\limits_{1\le k,j\le n}
           \frac{ (-\sigma+m_j-n+1)_{n-k}}
                  { (\tau+m_j-n+2)_{n-k}}
$$
Up to a permutation of rows, it is a determinant of
the form described in Lemma 1.2
with
$$x_j=m_j,\qquad a_j=-\sigma-n+j,\qquad \tau=-n+j+1
.$$
After a simple rearangement of the factors,
we obtain the required result
\hfill $\square$

\smallskip

{\bf 2.5. Hermitians forms defined by kernels.}
First, recall some standard facts on
characters of compact groups, for details
see, for instance, \cite{Kir}, 9.2, 11.1.

Let $K$ be a compact Lie group
equipped with the Haar measure $\mu$, we assume that
the measure of the whole group is 1.
Let
$\pi_1$, $\pi_2$, \dots
be the complete collection of pairwise distinct
irreducible representations.
Let $\chi_1$, $\chi_2$, \dots be their characters.
Recall the orthogonality relations
\begin{equation}
\langle \chi_k,\chi_l \rangle_{L^2}:=
\int_K \chi_k(h)\ov{\chi_l(h)}\,d\mu(h)=
\delta_{k,l}
\end{equation}
and
\begin{equation}
\chi_k*\chi_l= \begin{cases}
\frac 1{\dim\pi_k} \chi_k,&\qquad\text{if $k=1$}\\
                 0,&\qquad \text{if $k\ne l$}
\end{cases}
\end{equation}
where $*$ denotes the convolution on the group,
$u*v(g)=\int u(gh^{-1})v(h)d\mu(h)$.

We consider the action of the group $K\times K$
on $K$ by the left and right shifts
$(k_1,k_2):\,\,g\mapsto k_1gk_2$.
The representation $K\times K$
in $L^2(K)$ is a multiplicity free direct sum of
irreducible represntations having the form
$\pi_k\otimes\pi_k^*$, where $\pi_k^*$
is the dual representation,
\begin{equation}
L^2(K)\simeq \bigoplus_k \rho_k\otimes \rho^*_k
.\end{equation}
Each distribution $f$ on $K$ is a sum
of "elementary harmonics"
$$f=\sum\nolimits_k f^{(k)},\qquad f_k\in V_k.$$
The summands of this sum correspond to the decomposition
(2.14).

The projector to a subspace $\rho_k\otimes\rho_k^*$
 of $k$-th elementary harmonics
is the convolution with the corresponding character,
\begin{equation}
f^{\{k\}}=\frac 1{\dim \pi_k} f*\chi_k
\end{equation}
(in particular, $f^{(k)}$ is smooth).

\smallskip

 For  a central distribution $\Xi$ on $K$,
consider the Hermitian form
$$
\langle u,v\rangle=\iint_{K\times K} \Xi(h,g)\,u(h)\ov{v(g)}
\,d\mu(u)\,d\mu(g)
.$$
Consider the expansion of $\Xi$ in characters
$$\Xi=\sum_k c_k\chi_k.$$

{\sc Lemma 2.4.}
\begin{equation}
\langle u, v\rangle =
\sum_k\frac{c_k}{\dim\pi_k}
\int_{\U(n)} u^{\{k\}}(h)\ov{ v^{\{k\}}(h) }
\,d\mu(h)
.\end{equation}

\smallskip

{\sc Proof.} We can assume $u=\chi_k$, $v=\chi_l$.
We evaluate
$$
\iint_{K\times K} \sum_k c_k\Xi(gh^{-1})\,
     \chi_k(h) \ov\chi_l(g)\,d\mu(g)\,d\mu(h)
.$$
using (2.12) and (2.13).\hfill $\square$

\smallskip

{\bf 2.6. Positivity.}
Let $\Re(\sigma+\tau)<1$.
Consider the sesquilinear form on $C^\infty(\U(n))$
given by
\begin{equation}
\langle q, r\rangle_{\sigma,\tau}=
\iint_{\U(n)\times\U(n)} L_{\sigma,\tau}(g,h) q(g) \overline {r(h)}
\,dg\,dh
,\end{equation}
where the kernel $L_{\sigma,tau}$ is the same as above.
Obviously, for fixed $q$, $r$, this expression admits
a meromorphic continuation
in $\sigma$, $\tau$ to the whole $\C^2$.
Moreover, Theorem 2.3 allows to write an
explicit expression for this
continuation.
Expanding $q$ and $r$ in elementary  harmonics
$$
q(h)=\sum_\bm \kappa_\bm(h),\qquad r(h)=\sum_\bm \theta_\bm(h)
,$$
we obtain (see Lemma 2.4)
$$
\langle q, r\rangle_{\sigma,\tau}=
\sum_\bm \frac{c_\bm}{\dim \pi_\bm}
\int_{\U(n)}\kappa_\bm(h) \ov{ \theta_\bm(h)}\, d\mu(h)
,$$
where the meromorphic expressions for $c_\bm$ were obtained
in Theorem 2.3.

  If $\sigma$, $\tau\in \R$, then our kernel  $L_{\sigma,\tau}$
is Hermitian, i.e.,
$L_{\sigma,\tau}(h,g)=\ov{L_{\sigma,\tau}(g,h)}$,
or eguivalently
$$\langle q, r\rangle_{\sigma,\tau}=
\ov{\langle r, q\rangle}_{\sigma,\tau}
$$

\smallskip

{\sc Corollary 2.5.} {\it For $\sigma, \tau\in\R\setminus \Z$,
the inner product (2.17) is positive definite (up to a sign),
iff fractional parts of $-\sigma-n$ and $\tau$ are equal.}

\smallskip

The domain of positivity is the union
of the dotted squares on Figure 1.

For $\sigma$, $\tau$ satisfying this corollary, denote by
$H_{\sigma,\tau}$ the completion of $C^\infty(\U(n))$
with respect to our inner product.

\smallskip

{\bf  2.7. Action of $\U(n,n)$ on the space $\U(n)$.}
Consider the linear space $\C^n\oplus\C^n$ equipped
with the indefinite Hermitian form
$$
\{v\oplus w, v'\oplus w'\}=\langle v,v'\rangle_{\C^n\oplus 0}-
                             \langle w,w'\rangle_{0\oplus\C^n}
$$
where $\langle \cdot,\cdot\rangle$ is the standard inner product
in $\C^n$.
Denote by $\U(n,n)$ the group  of linear operators
in $\C^n\oplus\C^n$ preserving the  form $\{\cdot,\cdot\}$.
We write elements of this group as block
$(n+n)\times(n+n)$ matrices
$g:=\begin{pmatrix}
 \alpha&\beta\\ \gamma&\delta\end{pmatrix}$.
By definition, these matrices satisfy the condition
\begin{equation}
g\begin{pmatrix} 1&0\\0&-1\end{pmatrix}  g^*=
\begin{pmatrix} 1&0\\0&-1\end{pmatrix}
.\end{equation}

For $h\in\U(n)$, consider its graph $\graph(h)$ in $\C^n\oplus\C^n$.
It is an $n$-dimensional linear subspace, consisting of
all vectors $z\oplus zh$, where
a vector-row $z$ ranges in $\C^n$.
Since $h\in\U(n)$,
the subspace $\graph(h)$ is isotropic%
\footnote{A subspace $V$ in a linear space is
{\it isotropic} with respect to an Hermitian (or bilinear) form $Q$ if
$Q$ equals 0 on $V$.}
with respect to our Hermitian form $\{\cdot,\cdot\}$.
Conversely, any $n$-dimensional isotropic subspace
in $\C^n\oplus\C^n$ is a graph
of a unitary   operator $h\in\U(n)$.

Thus we have one-to-one correspondence between
the group $\U(n)$ and the Grassmannian of $n$-dimensional isotropic
subspaces in $\C^n\oplus\C^n$.

The group $\U(n,n)$ acts on the Grassmannian
in an obvious way, and hence $\U(n,n)$
acts on the space $\U(n)$.
An explicit formula for the latter action can be easily written:
\begin{equation}g=
\begin{pmatrix}\alpha&\beta\\ \gamma&\delta\end{pmatrix}:\quad
h\mapsto h^{[g]}:=(a+zc)^{-1}(b+zd)
,\qquad\qquad h\in\U(n),\,\,g\in\U(n,n)
\end{equation}

{\sc Lemma 2.6.} a) {\it For the Haar measure $\mu(h)$ on
$\U(n)$, we have}
$$
\mu\bigl(h^{[g]}\bigr)=
|\det\nolimits^{-2n}(\alpha+z\gamma)| \cdot \mu(h)
$$

b) {\it The kernel $L_{\sigma,\tau}$ satisfies the identity}
\begin{equation}
L_{\sigma,\tau} (u^{[g]},v^{[g]})
=L_{\sigma,\tau}(u,v)
\det(\alpha+u\gamma)^{\{\sigma|\tau\}}
\det(\alpha+v\gamma)^{\{\tau|\sigma\}}
.\end{equation}

\smallskip

{\sc Proof.}
a) The differential   of the map $u\mapsto u^{[g]}$
is given by
\begin{equation}
du\mapsto (a+u\gamma)^{-1}\, du\, (-\gamma u^{[g]} +\delta)
,\end{equation}
see, for instance, \cite{Ner-beta}, Lemma 1.1.

The rational map $u\mapsto u^{[g]}$ is defined on the space of all
complex $n\times n$ matrices. By (2.21),
its complex Jacobian
is
$$
J(g,u):=\det\nolimits^{-n}(a+zc)^{-n}
\det\nolimits^n  (-\gamma u^{[g]} +\delta)
.$$
Hence the real Jacobian on the space of all matrices
is  $|J(g,u)|^2$, and the Jacobian of the map
$\U(n)\to\U(n)$ is $|J(g,u)|$.
It can easily be checked (see \cite{Ner-beta}, Lemma 1.2)
 that
$$
\det (-\gamma u^{[g]} +\delta) =\det(\alpha+z\gamma)^{-1}\det(g)
.$$
By (2.18), $|\det g|=1$, and this finishes proof.

\smallskip

b) A direct calculation.

\smallskip

{\bf  2.8. Unitary representations of $\U(n,n)$.}
Let $\U(n,n)$  acts in the space of $C^\infty$-functions
on $\U(n)$
by the operators
\begin{equation}
\rho_{\sigma,\tau}(g) F(h)=
F(h^{[g]})\det\nolimits^{\{-n-\sigma|-n-\tau\}}(\alpha+h\gamma)
.\end{equation}

{\sc Remark.} Let us explain the sence
 of the complex power in this formula. It can easily be checked with (2.18),
that
$\|\gamma\alpha^{-1}\|<1$. Hence, for all matrices $h$ satisfying
$\|h\|\le 1$, the matrix $\alpha+h\gamma$ is invertible.
Hence the function
$$\ln\det(\alpha+h\gamma)
$$
 has a countable family of continuous
branches on the set $\|h\|\le 1$ and in particular on $\U(n)$.
We define
\begin{multline}
\det\nolimits^{\{-n-\sigma|-n-\tau\}}(\alpha+h\gamma)
 :=
\exp\bigl\{-(n+\sigma) \ln\det(\alpha+h\gamma)
   -\\ -(n+\tau)\ov{\ln\det(\alpha+h\gamma)}-2\pi i k(\sigma-\tau)\bigr\}
.
\end{multline}
Thus, we can think that for each $g\in\U(n,n)$
formula (2.22)   defines a countable family of operators
$\rho_{\sigma,\tau}(g)$, they differs one from
another by  constant factors
$\exp\{2\pi i k(\sigma-\tau)\}$.
These operators define a projective representation
(see \cite{Kir},14) of
the group $\U(n,n)$
$$
\rho_{\sigma,\tau}(g)\rho_{\sigma,\tau}(g')=
\lambda(g,g') \rho_{\sigma,\tau} (gg'),
\qquad \lambda(g,g')\in\C
.$$
Equivalently, we can consider the multi-valued
 operator-valued function
$\rho_{\sigma,\tau}(g)$  on $\U(n,n)$ as a single-valued
function on the universal covering group $\U(n,n)^\sim$
of $\U(n,n)$. Then $\rho_{\sigma,\tau}(g)$
became a linear representation of $\U(n,n)^\sim$.

  If $(\sigma-\tau)\in\Z$, then
$\exp\{2\pi i k(\sigma-\tau)\}$
and (2.23) is a well defined   single-valued
expression. In this case $\rho_{\sigma,\tau}$
is a linear representation of $\U(n,n)$.

\smallskip

{\sc Proposition 2.7.} {\it  The operators
$\rho_{\sigma,\tau}(g)$ preserve the form
$\langle\cdot,\cdot\rangle_{\sigma,\tau}$.}

\smallskip

{\sc Proof.} First, let $\Re(\sigma+\tau)<1$.
Substitute $h_1=u_1^{[g]}$, $h_2=u_2^{[g]}$ to the integral
$$
\iint_{\U(n)\times\U(n)}
L_{\sigma,\tau}(h_1,h_2)\,q(h_1)\ov{r(h_2)}\,d\mu(h_1)\,d\mu(h_2)
.$$
By Lemma 2.5, we obtain
\begin{multline*}
\iint_{\U(n)\times\U(n)}
L_{\sigma,\tau}(u_1,u_2)
\det(\alpha+u_1\gamma)^{\{-\sigma|-\tau\}}
|\det(\alpha+u_2\gamma)|^{-\tau|-\sigma\}}
\times\\ \times
\,q(h_1)\ov{r(h_2)}
                            |\det(\alpha+u_1\gamma)|^{-2n}
                            |\det(\alpha+u_2\gamma)|^{-2n}
\,d\mu(u_1)\,d\mu(u_2)
.\end{multline*}
Thus, our operators preserve the form
$\langle\cdot,\cdot\rangle_{\sigma,\tau}$.

 For general $\sigma$, $\tau$,
we consider the analytic continuation.

\smallskip

{\sc Corollary 2.8.} {\it For $\sigma$, $\tau$
satisfying the positivity conditions of Corollary 2.5,
the representation $\rho_{\sigma,\tau}$ is unitary.}

\smallskip

{\bf 2.9. Some remarks.}
a) {\it The case $n=1$}, $\sigma=\tau=s-1$
gives precisely the complementary series
of representations of $\SL(2,\R)\sim\SU(1,1)$
described above in 0.1.

\smallskip

b) The representations $\rho_{\sigma,\tau}$ are
very degenerated in the followng
sense. For any irreducible unitary
representation of a semisimple group $G$,
its restriction to the maximal compact subgroup
$K$ has a spectrum with finite multiplicities,
but usually these multiplicities are not bounded.

In our case, i.e., $G=\U(n,n)$, $K=\U(n)\times\U(n)$,
the restriction of $\rho_{\sigma,\tau}$ to $K$
is the multiplicity free sum $\rho_\bm\otimes\rho_\bm^*$
(thus, only few representations
of $K$ are present in the spectrum).

\smallskip

c) {\it Shifts of parameters.}
 For integer $k$, the {\it projective} representations
$\rho_{\sigma+k,\tau-k}$ and $\rho_{\sigma,\tau}$
are equivalent. The intertwining operator is
the multiplication by the determinant
$$
F(h)\mapsto F(h)\det(h)^k
.$$
This operator also defines an isometry
of forms
$L_{\sigma+k,\tau-k}$ and $L_{\sigma,\tau}$

\smallskip

d) {\it Symmetry.}  Representations $\rho_{-n/2-p,-n/2-q}$
and $\rho_{-n/2+p,-n/2+q}$   are dual.
The invariant pairing is given by the formula
\begin{equation}
(F_1,F_2) \mapsto \int_{\U(n)} F_1(h) F_2(h)\,d\mu(h)
.
\end{equation}

For verification of this statement, we substitute
$h\mapsto h^{[g]}$ and apply the formula for the Jacobian.

In particular, the point $(\sigma,\tau)=(n/2,n/2)$
corresponds to a representation of
a unitary principle series of $\U(n,n)$
(it is unitary in the space $L^2(\U(n))$).

%On Fig.1, for odd $n$, this point is inside square,
%and for even $n$ it is a joint vertex of two squares.

\kartinka

e) {\it Another symmetry.}
 The representation $\rho_{\tau,\sigma}$ is complex cojugate
to $\rho_{\sigma,\tau}$

\smallskip

f) {\it Problem of unitarisability of
subquotients.}
 Corollary 2.8 gives a classification of unitary representations
among $\rho_{\sigma,\tau}$.
But for integer $\sigma$ or integer $\tau$,
the representation $\rho_{\sigma,\tau}$
can contain a unitary subrepresentation,
 a unitary factor-representation, or a unitary sub-factor.

\smallskip

g) {\it Unitary highest weight representations
and Sahi's unipotent representations.}
The kernels
$$L_{\sigma,0}(u,v)=\det\nolimits^\sigma(1-uv^{-1})$$
are well-known, and they define highest weight representations.
By a well-known theorem
(Berezin, Gindikin, Rossi--Vergne, Wallach),
 the form $L_{\sigma,0}$ is a nonnegative
 Hermitian form
iff
\begin{equation}
\sigma<-n+1 \qquad \text{and for $\sigma=-n+1$, $-n+2$, \dots, 0}
.
\end{equation}

First, let  $\sigma<-(n-1)$,
If $m_n$ is negative, then the factor
$\Gamma(0+m_j+1)$ in the denominator
of (2.10) is infinity, and hence (2.10) is 0.
If $m_n\ge 0$, then all other coeffitients $c_\bm$
have the same sign (all signs are positive or all are negative).
It is easy to observe this from (2.10)
for noninteger $\sigma$. For integer $\sigma$,
this follows from the limit considerations:
numerator and denominator
in (2.19) have poles of the same order,
 the limit as $\sigma\to-l$
have to be nonnegative as a limit of a nonnegative function.

Second, let $\sigma$ be in the discrete part of the set
(2.25), $\sigma=-n+\theta$. Then the factor $\prod\Gamma(\alpha+j)$
has a pole of order $n-j$. The factor
$\prod\Gamma(\sigma+n-m_j)$ has a pole of order $\ge n-j$
(since $m_n\ge 0$, $m_{n-1}\ge 1$, \dots, $m_1\ge n-1$).
Hence the ratio is finite, it is nonvanishing if the order
of a pole of denominator is precisely $n-j$.
This happens iff $m_n=0$, $m_1=1$,\dots, $m_{n-j+1}=j-1$.
After this, it remains to follow signs in (2.11).

\smallskip

{\sc Remark.} This consideration almost coincides
with the original Berezin's proof \cite{Ber}.
Hua Loo Keng \cite{Hua} have obtained an expansion of the
kernel $L_{\sigma,0}$ in a series of characters
(it is a partial case] of our Theorem 2.3),
and Berezin checked signs in this expansion.

\smallskip

In integer points lying in the strip
$0\ge \sigma+\tau \ge -n+1$ also there are
located
unipotent representations
\cite {Sah4}, see below Section 5.

\smallskip

h) {\it Some other unitary subquotients.}
 Consider a representation $\rho_{\sigma,\tau}$
lying on the boundary of the dotted domain on Fig.1.
For definiteness, assume $\tau=0$.
Obviously, each term of its Jordan--Holder series
is unitarizable (since a limit of positive inner products
is positive).
Denote by $Y_j$ the set of all the signatures
$\bm$ satisfying the condition
$$m_{n-j}\ge 0,\qquad m_{n-j+1}<0$$
(i.e., precisely $j$ terms of the signature are negarive).
For $\bm\in Y_j$, the expression
$c_\bm$
(2.19) has a zero of order $j$.

Denote by $W_j$, the subspace in $C^\infty(\U(n))$ 
spanned by all harmonics with signatures 
lying in $\cup_{k\ge j} Y_j$. We obtain
an invariant filtration
$$C^\infty(\U(n))=W_0\supset W_1\supset\dots \supset W_n$$
The representation of $\U(n,n)$
in each sub-factor
$W_j/W_{j+1}$ is unitaty.

%Their restrictions to the subgroup
%$\U(n)\times\U(n)$ have  simple spectra.

%Apparently, these limit reprentations are
%related to Dvorsky and Sahi works \cite{DS}.

\smallskip

i) {\it Matrix Sobolev spaces of an arbitrary order.}
 Denote
$$s=-\sigma-\tau+n.$$

 Let $F$ be a distribution on $\U(n)$,
let $F=\sum F_\bm$ be its expansion
 in a series of elementary harmonics.
We have
\begin{multline}
F\in H_{\sigma,\tau}
\qquad \Longleftrightarrow\qquad
\sum_\bm \frac{c_\bm}{\dim \pi_\bm} \|F_\bm\|^2_{L^2}<\infty
\qquad \Longleftrightarrow
\\  \Longleftrightarrow
\qquad
\sum_\bm\Bigl\{  \|F_\bm\|^2_{L^2}
\prod_{j=1}^n (1+|m_j|)^s\Bigr\}<\infty
,\end{multline}
where  $\|F_\bm\|_{L^2}$
 denotes
$$
\|F_\bm\|_{L^2}:=
\Bigl(\int_{\U(n)} |F_\bm(h)|^2\,d\mu(h)\Bigr)^{1/2}
$$
Our Hermitian form   defines a norm only in the case
$|s|<1$, but (2.26) has sence for arbitrary real $s$,
and thus {\it we have a possibility to define
a Sobolev space ${\mathcal H}_s$ on $\U(n)$ of an arbitrary order}.

Author do not know applications
of this remark, but it seems that
 it can be useful in two  following
situations.

First, a reasonable harmonic
 analysis related to semisimple
Lie groups is the analysis
of unitary represenations.
 But near 1980 Molchanov observed that there are
many identities with special function that admits
interpretations  on "physical level of rigor" as
formulae of nonunitary  harmonic analysis.
Up to now, there are no reasonable interpretations
of this phenomenon (see, for instance
\cite{vDM}, see also
\cite{Ner-berezin}, Section 1-32 and formula (2.6)--(2.15) ).
In particular, we do not know
reasonable functional spaces that can be
place of action of this analysis.
 It seems that  our spaces
$H_s$ can be possible candidates.

Second, natural integral operators
 in the noncomutative harmonic  analysis
seem similar to pseudo-differential operators,
but they
are not pseudo-differential operators in the usual sence.
In particular, they are not well compatible with the
 standard  scales of functional spaces.
It can happened that our spaces $\cH_s$
can be more reasonable in this situation.

\smallskip

j) {\it An identity for formal series.}
We have  (see notation (0.4))
$$
\bigl|\sin (\phi/2)\bigr|^{\sigma+\tau}
e^{(\sigma-\tau)\phi}=
\frac{2^{\sigma+\tau}\Gamma(\sigma+\tau+1)}
{\Gamma(\sigma)\Gamma(\tau)}
\sum_{m=-\infty}^\infty
\frac{(-1)^m e^{im\phi}}
{(\tau+1)_m(\sigma+1)_{-m}}
.$$
Introduce  new variables $x_p=e^{i\phi_p}$.
We can rewrite Theorem 2.3 as the following
identity for formal series
\begin{multline*}
\prod_{p=1}^n
\Bigl\{\sum_{m=-\infty}^\infty
\frac{(-1)^m x_p^l}
{(\tau+1)_m(\sigma+1)_{-m}}
         \Bigr\}
=                \\=
\frac{\prod_{k=1}^{n-1} (\sigma+\tau+1)_k}
     {(\sigma)_n^n}
\cdot\sum_\bl \Bigl[\prod_{j=1}^n\frac{(-1)^{j_j}}
               {(\tau+1)_{l_j}(\sigma+n)_{-l_j}}
                 \Bigr]
              \frac{\bigl\{ \det_{1\le j,p\le n} x_p^{l_j}\bigr\} }
              {\bigl\{ \det_{1\le j,p\le n} x_p^{j-1}\bigr\} }
\end{multline*}

\bigskip

{\bf\large 3. Orthogonal groups}

\medskip

\stepcounter{sec}
\setcounter{equation}{0}

{\bf 3.1. Definition of the kernel.}
We consider the (disconnected) group $\OO(2n)$ as a basic object%
\footnote{It is also possible to concider its connectd subgroup
$\SO(2n)$;
in this case we must  replace the group $\OO(2n,2n)$ below by
it connected component $\SO_0(2n,2n)$, consider
a connected component of the Grassmannian and also do some
obvious minor changes.}.
Each element of this group can be reduced by
a conjugation $g\mapsto hgh^{-1}$
to the block diagonal form with $2\times 2$-blocks
\begin{equation}
\!\!\!\!
\begin{pmatrix}
A(\phi_1)&0&\dots&0\\
0&A(\phi_2)&\dots&0\\
\vdots&\vdots&\ddots&\vdots\\
0&0&\dots&A(\phi_n)
\end{pmatrix}
\text{, where
$A(\phi_j)=
\begin{pmatrix}
\cos\phi_j &\sin \phi_j\\
-\sin\phi_j&\cos\phi_j
\end{pmatrix}$}
.\end{equation}

The collection $(\phi_1,\dots,\phi_n)$ is uniquely determined
by $g$ modulo permutations and arbitrary transformations
$\phi_j\mapsto -\phi_j$.
The numbers $e^{\pm i\phi_j}$ coincide with the eigenvalues
of $g$.

We define the function $\ell_\lambda(g)$
on $\OO(2n)$
by
$$
\ell_\lambda(g)=\begin{cases}
                           \det\bigl((1-g)/2\bigr)^\lambda,
                                  &, \text{if $g\in\SO(2n)$;}\\
                          0,& \text{if $g\notin\SO(2n)$.}
               \end{cases}
$$
the derminant is nonegative
 and hence its complex powers are well-defined.
In the terms of the eigenvalues $e^{\pm i\phi_j}$,
the function $\ell_\lambda(g)$
coincides with
$$
\ell_\lambda(g)=\prod_{j=1}^n \bigl|\sin(\psi_j/2)\bigr|^{2\lambda}
.$$

We define the kernel $L_\lambda(\cdot,\cdot)$ on $\OO(2n)$ as
$$
L_\lambda(g,h)=\ell_\lambda(gh^{-1})
.$$

{\bf 3.2. Characters,} see \cite{Wey}.
Irreducible representations $\pi_\lambda$ of $\SO(2n)$
are parametrized  by collections
of numbers
$$
(\bl,\epsilon):\quad l_1>l_2>\dots>l_n> 0, \qquad \epsilon=\pm 1,
$$
or
$$
\bl:\quad l_1>l_2>\dots>l_n= 0.
$$
Formulae for the characters
 are slightly different in these two cases.

If $l_n=0$, we have
$$
\chi_\bl(g)=\frac{\det_{1\le k,m\le n}
        \bigl\{ \cos(l_k\phi_m)\bigr\}}
            { \det_{1\le k,m\le n}\bigl\{ \cos(k-1)\phi_m\bigr\}}
.$$

For $l_n\ne 0$, we have
$$
\chi_\bl^\epsilon(g)=\frac{\det_{1\le k,m\le n}
        \bigl\{ \cos(l_k\phi_m)\bigr\}+
           \epsilon\det_{1\le k,m\le n}\bigl\{ \sin(l_k\phi_m)\bigr\} }
            {2 \det_{1\le k,m\le n}\bigl\{ \cos(k-1)\phi_m\bigr\}}
.$$

 Denote by $J$ the diagonal matrix $\in\OO(2n)$
having $(2n-1)$ entries   $(+1)$
and one $(-1)$. The map $h\mapsto JhJ^{-1}$ is an interior
automorphism of $\OO(2n)$ and an exterior automorphism
of $\SO(2n)$.  The representations $\rho^\pm_\bl$
of $\SO(2n)$
corresponding to the characters $\chi^\pm_\bl$ are twins
in the following sence
$$
\pi^\pm_\bl(JhJ^{-1})=\pi^\mp_\bl(h)
.$$
Also the substitution
$$
\phi_1\mapsto \phi_1,\qquad \dots\qquad\phi_{n-1}\mapsto \phi_{n-1},
\qquad
\phi_n\mapsto-\phi_n
$$
changes $\chi_\bl^+$ and  $\chi_\bl^-$.

For $l_n=0$, the representation $\pi_\bl$ is
stable with respect to the exterior
automorphism $h\mapsto JhJ^{-1}$.

This digression about two types of characters
 is almost non-essential for us, since really we
will consider arbitrary
$$\bl:\quad l_1>l_2>\dots>l_n\ge 0, \qquad  $$
and the
 functions $\chi_\bl$    given by
\begin{equation}
\!\!
\chi_\bl(g):=
\frac{\det_{1\le k,m\le n}
        \bigl\{ \cos(l_k\phi_m)\bigr\}}
            { \det_{1\le k,m\le n}\bigl\{ \cos(k-1)\phi_m\bigr\} }
=\begin{cases}\chi_\bl(g),&\,\text{if $l_n=0$,}
\\
\chi_\bl^+(g)+\chi_\bl^-(g), &\,\text{if $l_n>0$.}
\end{cases}
\end{equation}

{\sc Remark.}
The functions (3.2) are precisely  restrictions
of the characters of the disconnected group $\OO(2n)$
to the connected group $\SO(2n)$.

\smallskip

Let $F$ be a central functions on $\SO(2n)$.
Let $f$ be its restriction
to the maximal torus, i.e., to the set of matrices
having the form (3.1).
The {\it Weyl integration formula} for central functions
on $\SO(2n)$
has the form
\begin{multline*}
\int\limits_{\SO(2n)} F(g)\/d\mu(g)=\\=
\frac 1{\pi^n n!}
\int\limits_{0<\psi_1<2\pi} \dots\int\limits_{0<\psi_n<2\pi}
f(\phi_1,\dots,\phi_n)
\Bigl(
    \det_{1\le k,m\le n}\bigl\{ \cos(k-1)\phi_m \bigr\}
\Bigr)^2\, d\phi_1\dots d\phi_n
.\end{multline*}

{\bf 3.3. Expansion of $\ell_\lambda$ in characters.}

\smallskip

{\sc Theorem 3.2.} {\it Let $\lambda< 1/2$.
For $g\in\SO(2n)$,
$$
\ell_\lambda(g)=\sum_{\bl:\, l_n=0} c_\bl \chi_\bl
+\frac12  \sum_{\bl:\, l_n>0} c_\bl (\chi^+_\bl+\chi_\bl^-)
,$$
where}
\begin{multline}
c_\bl=
(-1)^{n(n-1)/2}
2^{2n\lambda+1} \pi^{-n} %\Gamma^n(2\lambda+1)
\,\sin^n \pi \lambda
%\cdot\prod_{k=1}^{2n-1}(2\lambda+k)^{[n-k/2]}
\prod_{k=1}^n \Gamma(2\lambda+2k-1)            %\nonumber
\times \\ \quad \times
\prod_{1\le\alpha<\beta\le n}(l_\alpha^2-l_\beta^2)\cdot
\prod_{j=1}^{n}
   \frac{ \Gamma(l_j-\lambda-n+1)}
        {\Gamma(l_j+\lambda+n)}
=
\end{multline}
\begin{multline}
=
(-1)^{n(n-1)/2}
2^{2n\lambda+1} %\Gamma^n(2\lambda+1)\,
\prod_{k=1}^n \Gamma(2\lambda+2k-1)
%\cdot\prod_{k=1}^{2n-1}(2\lambda+k)^{[n-k/2]}%\nonumber
\times \\ \quad \times
\frac
{(-1)^{\sum l_j}\prod_{1\le\alpha<\beta\le n}(l_\alpha^2-l_\beta^2)}
{\prod_{j=1}^{n}\Gamma(-l_j+\lambda+n)\Gamma(l_j+\lambda+n))
}
.\end{multline}

{\sc Proof.}   First,
$\ell_\lambda(JhJ^{-1})=\ell_\lambda(h)$. Hence
$$
\langle \ell_\lambda,\chi^+_\bl\rangle_{L^2}
= \langle \ell_\lambda,\chi^-_\bl\rangle_{L^2}
=\frac12  \langle \ell_\lambda,\chi_\bl\rangle_{L^2}
$$
Evaluating the last expression, we obtain
\begin{multline*}
\frac 1{\pi^n n!}
\int\limits_{0<\phi_1<2\pi}
\dots
\int\limits_{0<\phi_n<2\pi}
\prod_{k=1}^n \bigl|\sin(\phi_k/2)\bigr|^{2\lambda}
\times \\ \times
\det\limits_{1\le j,k\le n}\bigl\{ \cos l_j \phi_k\bigr\}
\det\limits_{1\le j,k\le n}\bigl\{\cos (j-1)
\phi_k\bigr\}
\prod_{k=1}^n d\phi_k
.\end{multline*}
Applying Lemma 2.2, we transform this integral to
$$
\pi^{-n}\det I(j,m),
$$
where
$$
I(j,m)=\int_0^{2\pi}\bigl|\sin(\phi/2)\bigr|^{2\lambda}
 \cos (l_j\phi)
\,  \cos [(m-1)\phi] \,d\phi
.$$
Expanding the product of cosines, we obtain a sum of two
integrals of the form (2.9)
$$
\int_0^\pi \bigl|\sin \phi\bigr|^{2\lambda}
   \Bigl[\cos 2(l_j+m-1)\phi  +   \cos 2(l_j-m+1) \phi\Bigr]\,d\phi
=
$$
\begin{multline*}
=
\pi 2^{-2\lambda}
 \Gamma(2\lambda+1)
\Bigr[\frac{(-1)^{l_j+m-1} }
  {\Gamma(\lambda+l_j+m)\Gamma(\lambda-l_j-m+2)}
+\\+
   \frac{(-1)^{l_j+m-1}     }
  {\Gamma(\lambda+l_j-m+2)\Gamma(\lambda-l_j+m)}
\Bigr]
.\end{multline*}
After simple transformations, we get
\begin{multline*}
2^{-2\lambda} \Gamma(2\lambda+1)\sin\lambda\pi
\Bigl[\frac{\Gamma(-1-\lambda+l_j+m)}
           {\Gamma(\lambda+l_j+m)}
     +
      \frac{\Gamma(1-\lambda+l_j-m)}
           {\Gamma(2+\lambda+l_j-m)}
\Bigr]         =
\\ =
\frac{\Gamma(2\lambda+1)\sin(\lambda\pi)}
{2^{2\lambda}}
\cdot\frac{\Gamma(l_j-\lambda)}
          {\Gamma(l_j+\lambda+1)}
\boxed{
\Bigl( \frac{(l_j-\lambda)_{m-1}}
             {(l_j+\lambda+1)_{m-1}}
+
\frac {(l_j+\lambda-m+2)_{m-1}}
      {(l_j-\lambda-m+1)_{m-1}}
\Bigr) }
\end{multline*}
The factors outside the boxed equation do not depend on $m$.
Thus, we must evaluate
$$
\det_{1\le j, m\le n}
\Bigl[ \frac{(l_j-\lambda)_{m-1}}
             {(l_j+\lambda+1)_{m-1}}
+
\frac {(l_j+\lambda-m+2)_{m-1}}
      {(l_j-\lambda-m+1)_{m-1}}
\Bigr]
.$$
A matrix element can be represented in the from
\begin{multline*}
h_{jm}=\frac{  (l_j-\lambda)(l_j-\lambda+1)\dots(l_j-\lambda+m-2)}
     {(l_j+\lambda+1) (l_j+\lambda+2)\dots (l_j+\lambda+m-1)}
+\\+
\frac{  (l_j+\lambda)(l_j+\lambda-1)\dots(l_j+\lambda-m+2)}
     {(l_j-\lambda-1) (l_j-\lambda-2)\dots (l_j-\lambda-m+1)}
.\end{multline*}
(in each numerator and denominator we have $(m-1)$ factor;
in particular for $m=1$ we have $h_{j1}=2$).
Thus we obtain a determinat evaluated in Lemma 1.4
with
$$
x_j=l_j,\qquad a_j=-\lambda+j-1,\qquad b_j=\lambda_j+j
.$$
After some rearangement  of the factors we
obtain the required result.

\smallskip

{\bf 3.4. Representations.} Now we can repreate
the considerations of 2.6-2.8.
Consider the linear space $\R^{2n}\oplus\R^{2n}$ equipped
with the indefinite bilinear form
$$
\{v\oplus w, v'\oplus w'\}=\langle v,v'\rangle_{\R^{2n}\oplus 0}-
                             \langle w,w'\rangle_{0\oplus\R^{2n}}
,$$
where $\langle \cdot,\cdot\rangle$ is the standard inner product
in $\R^{2n}$.
Denote by $\OO(2n,2n)$ the group  of linear operators
in $\R^{2n} \oplus\R^{2n}$ preserving this  form.

For $h\in\OO(2n)$ consider its graph $\graph(h)$ in
$\R^{2n}\oplus\R^{2n}$.
Since $h\in\OO(2n)$,
the subspace $\graph(h)$ is isotropic
with respect to our bilinear form $\{\cdot,\cdot\}$.

Thus we have one-to-one correspondence between
the group $\OO(2n)$ and the Grassmannian of $2n$-dimensional isotropic
subspaces in $\R^{2n}\oplus\R^{2n}$.

The group $\OO(2n,2n)$ acts on the Grassmannian
 and hence $\OO(2n,2n)$
acts on the space $\OO(2n)$.
The explicit formula for the latter action
is the same as above,
see (2.19).

We consider the inner product in $C^\infty(\OO(2n))$
given by
\begin{equation}
\langle F_1,F_2\rangle_\lambda=
\iint_{\OO(2n)\times\OO(2n)}
L_\lambda(h_1,h_2)\,F(h_1)\, \ov{F(h_2)}\,d\mu(h_1)\, d\mu(h_2)
.\end{equation}
This integral is convergent if $\lambda>-1/2$, for general
$\lambda$ we consider the analytic continuation.

\smallskip

{\sc Proposition 3.2.} a)
{\it The inner product  (3.5) is invariant with respect to
the transformations
$$
\rho_\lambda(g) F(h)=
F(h^{[g]})\det(\alpha+h\gamma)^{-2n+1-\lambda}
,$$
where $g\in\OO(2n,2n)$, $h\in\OO(2n)$.}

b)
{\it If $-(n-1)>\lambda >-n$, then
the inner product $\langle\cdot,\cdot\rangle_\lambda$
is positive definite. In other words, the representation
$\rho_\lambda$ in this case is unitary.}

\smallskip

{ \sc Proof.} Statement b) follows from Theorem 3.1.
For evaluation of the Jacobian,
see, for instance, \cite{Ner-beta}, 1.2--1.3.

\medskip

{\bf\large 4. Symplectic groups}

    \stepcounter{sec}
\setcounter{equation}{0}

\medskip

{\bf 4.1. Quaternionic matrices.}
We denote by $\H$ the quaternionic field.
Operators in the  quaternionic coordinate
 space $\H^n$ can written in the form
$$
v\mapsto v Q
,$$
where $Q$ is an $n\times n$-matrix with quaternionic elements,
and $v\in\H^n$ is a vector-row.
More formaly, these transformations are endomorphisms
of a left $n$-dimensional module over the field $\H$.

Let $Q$ be a quaternionic matrix. It defines the operator
$\H^n\to\H^n$, identifying $\H^n$ with $\R^{4n}$,
Hence we obtain an operator $Q_\R:\R^{4n}\to \R^{4n}$.
We define the determinant of $Q$ by
$$
\det Q=\sqrt[4]\det Q_\R.
$$
($\det Q_\R$ is a nonnegative real number,
hence $\det Q$ also is a positive real number).

A quaterninic matrix $Q$ is unitary, if it preserves
the Hermitian form
$$
\{v,w\}=\sum v_j \ov w_j
.$$
In other words
$$
QQ^*=1,
$$
where $Q^*$ is a matrix obtained from $Q$
by transposition and element-wise quaternionic conjugation;
if $q_{kl}$ are the matrix elements of $Q$, then the matrix elements
of $Q^*$ are $\ov q_{lk}$.

\smallskip

{\bf 4.2. Symplectic group,}
see \cite{Wey}.
The group of all $n\times n$ quaterninic unitary
matrices is called "symplectic group" or
"compact symplectic  group",
the standard notation is $\Sp(n)$.

Using a group conjugation $h\mapsto  g^{-1}hg$,
each element of $\Sp(n)$ can be reduced to
a diagonal matrix with the entries
$e^{ik\phi}$.
The collection $\psi_j$ is defined uniquely
 modulo permutations and change of signs:
$\psi_j\mapsto (-1)^{\sigma_j} \psi_j$.

Let $F$ be a central function on $\Sp(n)$,
let
$f$ be its restriction to the set of diagonal matrices.
The following {\it Weyl integration formula} holds
\begin{multline*}
\int_{\Sp(n)} F(g)\,d\mu(g)=
\frac1{\pi^n n!}
\int\limits_{0<\psi_1<2\pi}
\dots
\int\limits_{0<\psi_n<2\pi}
f(\phi_1,\dots,\phi_n)
\times \\ \times
\Bigl|\det_{1\le j,k \le n}\bigl\{ \sin (k\psi_j)\bigr\}\Bigr|^2
\prod_{j=1}^n d\phi_j
.
\end{multline*}

Finite-dimensional irreducible representations
of $\Sp(2n)$ are enumerated by collections
$\bl$ of integers having the form
$$
\bl:\quad l_1>l_2>\dots>l_n>0
.$$
The corresponding character is
$$
\chi_l(g)=
\frac {\det_{1\le j,k \le n} \bigl\{ \sin (l_k \psi_j)\bigr\} }
      {\det_{1\le j,k \le n} \bigl\{ \sin (k\psi_j)\bigr\}}
.$$

We define the function $\ell_\lambda$ on $\Sp(n)$
by
$$
\ell_\lambda(g)=\det\bigl((1-g)/2\bigr)^{2\lambda}=
\prod_{j=1}^n |\sin(\psi_j/2)|^{2\lambda}
.$$

\smallskip

{\bf 4.3. Expansion of $l_\lambda$.}

\smallskip

{\sc Theorem 4.1.} {\it For $\Re\lambda<1/2$,
$$
\ell_\lambda(g)=
\sum_\bl c_\bl \chi_\bl
,$$
where}
\begin{multline}
c_\bl=\pi^{-n} 2^{-2n\lambda} \sin^n \lambda\pi
\prod_{k=1}^n \Gamma(2\lambda+2k)
\times \\ \times
\prod_{\alpha=1}^n (2l_\alpha)
\prod_{0\le\alpha<\beta\le n} (l_\alpha^2-l_\beta^2)
\prod_{j=1}^n \frac{\Gamma(l_j-\lambda-k)}
          {\Gamma(l_j+\lambda+1+k)}
.\end{multline}

{\sc Proof} is very similar to the orthogonal case.
We must evaluate the determinat whose matrix elements
are
$$
I(k,j):=
\int_0^{2\pi}\bigl|\sin(\phi/2)\bigr|^{2\lambda}
\sin(l_j\phi)\, \sin(k\phi)\,d\phi
.$$
Expanding the product of sines, we obtain a sum of two
integrals of the form (2.9)
\begin{multline*}
\int_0^{\pi}\bigl|\sin(\phi/2)\bigr|^{2\lambda}
\bigl[\cos 2(l_j-k)\phi -\cos 2(l_j+k)\phi\bigr]
\,d\phi
=\\=
\pi 2^{-2\lambda}\Gamma(2\lambda+1)
\times\\ \times
\Bigl[\frac{(-1)^{l_j-k}}
   {\Gamma(l_j+\lambda-k+1)\Gamma(-l_j+\lambda+k+1)}-
      \frac{(-1)^{l_j+k}}
   {\Gamma(l_j+\lambda+k+1)\Gamma(-l_j+\lambda-k+1)}
\Bigr]
=\\=
\pi 2^{-2\lambda}\Gamma(2\lambda+1)
\Bigl[\frac{\Gamma(l_j-\lambda-k)}
   {\Gamma(l_j+\lambda-k+1)}+
      \frac{\Gamma(l_j-\lambda+k)}
   {\Gamma(l_j+\lambda+k+1)}
\Bigr]
.\end{multline*}
We transform the expression in the brackets to
$$
\frac{\Gamma(l_j-\lambda)}{\Gamma(l_j+\lambda)}
\cdot
\Bigl\{\frac{(l_j+\lambda-k+1)_k}
            {(l_j-\lambda-k)_k}
  - \frac{(l_j-\lambda)_k}
         {(l_j+\lambda+1)_k}
\Bigr\}
.$$
The factor in the front of brackets does not depend on $k$
and hence it is suffitient
to evaluate the determinant of the matrix
$$u_{kj}=\frac{(l_j+\lambda-k+1)_k}
            {(l_j-\lambda-k)_k}
  - \frac{(l_j-\lambda)_k}
         {(l_j+\lambda)_k}
.$$
This matrix has the form described in Lemma 1.3,
with
$$
x_j=l_j,\qquad a_k=\lambda+1-k,\qquad
b_k=-\lambda-k
.$$

{\bf 4.4. Unitary representations of $Sp(n,n)$.}
The pseudounitary quaternionic group $\Sp(n,n)$
is the group of quaternionic $(n+n)\times(n+n)$-matrices
preserving the Hermitian form
$$
\{ v,w\}:=
\sum_{j=1}^n v_j\ov w_j-\sum_{j=n+1}^{2n} v_j \ov w_j
.$$

We consider its action on $\Sp(n)$ by linear fractional
transformations (2.19) as above.
We also consider the representation $\rho_\lambda(g)$
of $\Sp(n,n)$ in the space $C^\infty(\Sp(n))$
given by
$$
\rho_\lambda(g) F(h)=
F(h^{[g]})\det(\alpha+h\gamma)^{-2n-1-\lambda}
$$

{\sc Proposition 4.2.}
a) {\it The operators $\rho_\lambda(g)$ preserve
the Hermitian form}
$$
\langle F_1, F_2\rangle_\lambda=
\iint_{\Sp(n)\times\Sp(n)}
\ell_\lambda(gh^{-1}) F_1(g)\ov{F(h)}\,
d\mu(g)\,d\mu(h)
.$$

b) {\it If $-n>\lambda>-n-1$, then the Hermitian form
$
\langle \cdot,\cdot\rangle_\lambda
$
is positive definite and thus
 the representation $\rho_\lambda$
is unitary.}

\smallskip

The nontrivial part of the statement is positivity
of the Hermitian form.
This follows from Theorem 4.1.

\medskip

{\bf\large 5. Unipotent representations.}

    \stepcounter{sec}
\setcounter{equation}{0}

\medskip

Here we discuss models of "unipotent" representations
of Sahi \cite{Sah4} and Dvorsky--Sahi, 
\cite{DS1}--\cite{DS2}.

\smallskip

{\bf 5.1. The case $\OO(2n,2n)$.}
In notation of Section 3,
we suppose that
\begin{equation}
\lambda=-n+\alpha\qquad
\text{$\alpha=1$, 2, \dots, $n$.}
\end{equation}
 
The first row in (3.3)
has zero of order $\alpha$
at our $\lambda$. The second row
in (3.3) has a pole of order
$\le \alpha$. Hence the total expression is nonzero
iff the order of the pole is precisely $\alpha$, i.e.,
$$
l_n=0,\quad, l_{n-1}=1,\quad\dots l_{n-\alpha+1}=\alpha-1
$$
Under this condition, all the coefficients
$c_\bl$ are positive.

Thus, for $\lambda$ having the form (5.1),
the inner product (3.5) is  non-negative definite
(and degenerated)
and the operators (3.6) are unitary with respect
to this inner product. 

\smallskip

{\sc Remark.} For $\alpha=0$, our representation
is the one-dimensional representation.
 
For $\alpha=1$, the representation
obtained in this way is an element of
the Molchanov degenerated discrete series,
see \cite{Mol},
now there exists a wide literature 
devoted to this representation, see,
for instance \cite{Kob}.

\smallskip

{\bf 5.2. Groups $\U(n,n)$.}
Now we assume 
\begin{equation}
\tau=0, \qquad \sigma=-n+\alpha,\qquad
\text{where $ \alpha=1$, \dots, $n-1$.}
\end{equation}
The coefficient
$c_\bm$ given by (2.10) is non-zero, if 
$\bm$ is contained in the union of 
the following disjoint sets $Z_j$,
$j=0$, 1, \dots, $n-\alpha$,
$$
Z_\theta: \bm\quad\text{has form}\quad
 (m_1,\dots, m_{n-\alpha-\theta}, \alpha-1,\alpha-2,\dots,0,
m_{n-\theta+1},\dots m_n)
$$

%Denote by $Tail$ the set of signatures 
%that do not contained in the union of $\cup  Z_j$.

Denote by $V_\bm$ be the $\U(n)\times\U(n)$-submodule
in $C^\infty(\U(n))$ corresponding a signature $\bm$,
see 2.5.

\smallskip

{\sc Proposition 5.1.} {\it The subspace
$$
W_{tail}:=\oplus V_\mu\subset C^\infty(\U(n)),\quad
\text{where $\mu\notin \cup Z_j$,}
$$
is a $\U(n,n)$-invariant subspace.}

\smallskip

b) {\it The quotient $C^\infty(\U(n)) /W_{tail}$
is a sum $n-\alpha+1$ submodules 
$$W_j=\oplus_{\mu\in Z_j} V_ \mu$$
The representation of $\U(n,n)$ in each
$V_j$ is unitary.} 
  
\smallskip

A proof is given below.

%Denote by $W_\theta$ a space of all $f\in C^\infty(\U(n))$,
%whose expansions in harmonics contain only
%terms numerated by
% $Z_\theta$.
% By $C^\infty_{tail}$ we denote the set of
% all $f\in C^\infty(\U(n))$ whose expansions
% in harmonic do not contain terms numerated by $\cup Z_j$.
 
% For a function $f^\in C^infty\U(n)$, we denote
% by $f_{[\theta]}$ its projection
% to  $W_\theta$.
 
% \smallskip

%{\sc Proposition 5.1.}
%a) {\it Subspace $W_{tail}$ is $\U(n,n)$-invariant.}%
%
%b) {\it Let us identify $C^\infty(\U(n))/W_{tail}$
%with $\oplus W_\theta$ in the natural way.
%In this sence, the subspace
% $W_\theta\subset C^\infty/W_{tail}$ is $\U(n,n)$-invariant,
%and 
%  the representation of $U(n,n)$
%in each $W_\theta$ is unitary}.

\smallskip

{\bf 5.3. Blow-up.}
The distribution $\ell_{\sigma,\tau}$
depends meromorphically in 
the two complex variables $\sigma$, $\tau$. 
Its poles and zeros are located at $\sigma\in\Z$
and in $\tau\in\Z$ and hence values of this distributions
at points $(\sigma,\tau)\in\Z^2$ generally
 are not uniquely defined.
 Passing to this point from different directions,
 we can obtain different limits.
 
It is convenient to formulate this more 
formally.

We fix a point $(\sigma,\tau)=(-n+\alpha,0)$
and introduce the new local coordinates near this point%
\footnote{This construction is the {\it blow-up}
of the plane $\C^2$ at the point $(-n+\alpha,0)$.}
by
\begin{equation}
\sigma=-n+\alpha+s\epsilon,\quad \tau=t\epsilon
\end{equation}
The new coordinates are defined up to the equivalence
\begin{equation}
\qquad (s,t,\epsilon)\sim (su,tu,\epsilon/u)
\end{equation}
We also think that 
\begin{equation}
|s|,\quad |t|,\quad |\epsilon|
\text{ are
suffitiently small,
and $(s,t)\ne (0,0)$.}
\end{equation}

If we replace $\epsilon=1/R$, then the 
collections $(s,t,R)$ are defined
up to the equivalence
$$
\qquad (s,t,\epsilon)\sim (su,tu,Ru)
$$
i.e., $s,t,u$ is a point of the projective plane.
Thus the set (5.5) can be consired as a subset in projective plane.
In the new coordinates, the point $(-n+\alpha,0)\in\C^2$
corresponds to the whole complex projective line 
$(s,t,0)$.

Thus replacing a neibourhood of  $(-n+\alpha,0)\in\C^2$
by the set (5.5), we obtain a new complex manifold,
denote it by $\wt\C^2$.

\smallskip

{\sc Proposition 5.2.}
{\it The distribution 
$$\ell^{s,t,\epsilon}(z):=\ell_{\sigma,\tau}(z)$$
 is a meromorphic distribution-valued function
 on $\wt\C^2$ in the domain (5.5).
 The unique pole (or order $n$) near $(s,t)=(0,0)$
  is the line
 $s+t=0$.}
 
\smallskip

{\sc Proof.}  Each Fourier coefficient $c_\bm$ in the formula
(2.10)
has the form
\begin{equation}
\epsilon^k 
\frac{t^j s^{n-\alpha-j}}
{(s+t)^{n-\alpha}}
 R_\bm(s,t,\epsilon)
\end{equation}
where the last factor is holomorphic and
nonvanishing near the line $(s,t,0)$,
and the powers $k$, $j$, $n-\alpha-j$ are nonnegative.

Thus, the Fourier coefficients $c_\bm$ of the distribution
$\ell^{s,t,\epsilon}$ are meromorphic in the our region.

We must verify that the derivatives 
$\partial \ell^{s,t,\epsilon}/\partial t $,
$\partial \ell^{s,t,\epsilon}/\partial s$, 
$\partial \ell^{s,t,\epsilon}/\partial \epsilon$
are well-defined distributions.

A central function on $\U(n)$ is a distribution, if
its Forier coefficient have at most polynomial growth
in $\bm$.
For (5.6), this can be verified by a direct tracing.

\smallskip

{\bf 5.4. The family of invariant Hermitian forms.
Proof of Proposition 5.1.}
Thus, for $(\sigma, \tau)=(-n+\alpha,0)$, 
we have the following family of Hermitian forms invariant
with respect to
the operators (2.22)
$$
\Lambda_{s,t}(f,g)=
\iint_{\U(n)\times\U(n)} 
\ell^{s,t,0}(zu^*)f(z)\ov {f(u)}\,d\mu(z)\,d\mu(u)
$$

Now we emphasis some additional properties
of the formula (5.6).

First, the exponent $k$ in (5.6) is zero iff $\bm\in\cup Z_j$.
Hence $W_{tail}$ is the common kernel of all the forms
$\Lambda_{s,t}$. This implies 
the statement a) of Proposition 5.1.

By the invariance, we subspaces $W_j$ are pairwise  orthogonal
with respect to all the forms $\Lambda^{s,t,0}$.

Secondly, $R_\bm(s,t,0)$ is a constant which depend only on
$\bm$ (since nonconstant summands in the linear factors of $R_\bm$
have the form $s\epsilon)$, $t\epsilon$).

Thirdly,
if $\bm\in Z_j$, then the exponent $j$ in (5.6) is
our $j$. This means, that the restriction
of $\Lambda^{s,t,0}$ to $W_j$
has the the form
\begin{equation}
\frac{t^j s^{n-\alpha-j}}
{(s+t)^{n-\alpha}} \Xi_j(f,g)
\end{equation}
where $\Xi(f,g)$ is independent of $t$, $s$.

The Fourier coefficients of $\Xi_j$
are the factors $R_\bm(s,t,0)$ from (5.6), 
and they can easily 
be written. It is easy to observe that 
all these coefficients have the same sign.

This finishes proof of Proposition 5.1.b).

\smallskip

{\bf 5.5. Expression for the kernel.} We can concider
the form $\Xi_j$ as a form on $C^ \infty(\U(n))$
extending it  to other harmonics as zero.

\smallskip

{\sc Proposition 5.3.} {\it The form $\Xi_j$
is defined  by the distribution} 
\begin{multline*}
\frac1 {n!}
\frac{\partial^j}{\partial t^j} 
(1+t)^{n-\alpha} \ell^{1,t,0}(zu^*)
\Bigr|_{t=0}=\\=
\frac1 {n!}
\lim_{\epsilon\to 0} \frac{\partial^j}{\partial t^j}
(1+t)^{n-\alpha}
\det(1-zu^*)^{\{\-n+\alpha+\epsilon|t\epsilon\}}
\Bigr|_{t=0,\epsilon=0}
\end{multline*}

{\sc Proof.} The first expression is obtained by differention
of the block decomposition (5.7)
of the Hermitian forms
$\Lambda_{s,t}$;
the right-hand side is a result of changing of the limit passing and
differentiation.   
  
 \smallskip 
 
{\sc Remark.} The sum (5.8)
contains the well-defined term
$\det(1-zu^*)^{-n+\alpha}:=\ell^{1,0,0}(zu^*)$
(that looks like a function)
and sum of distributions supported by
the submanifold $\mathrm{rk}(1-zu^*)=n-1$, $n-2$,
\dots, $n-\alpha$
(the analysis of the determinantal distributions
was done by Ricci, Stein \cite{RS}).

{\sf
Math. Phys.Group, Institute for the Theoretical and Experimental Physics,

Bolshaya Cheremushkinskaya, 25,

Moscow 117 259, Russia
\&
University of Vienna, Math. Dept.,

Nordbergstarsse, 15,
Vienna,  Austria
}

\tt neretin@mccme.ru

\end{document}